\documentclass[twoside,leqno]{article}

\usepackage[letterpaper]{geometry}

\usepackage{siamproceedings}

\usepackage[T1]{fontenc}
\usepackage{amsfonts}
\usepackage{graphicx}
\usepackage{epstopdf}
\usepackage{enumitem}
\usepackage{algorithmic}

\usepackage{tikz}
\usetikzlibrary{calc}
\usepackage{xcolor}
\usepackage[utf8]{inputenc}
\usepackage{amsmath,amssymb, bbm}
\usepackage{hyperref}
\hypersetup{
  colorlinks,
  linkcolor={red!50!black},
  citecolor={blue!50!black},
  urlcolor={blue!80!black},
  hypertexnames=false
}
\usepackage{cleveref}
\usepackage{nicematrix}
\usepackage{xspace}
\usepackage{thmtools}
\usepackage{thm-restate}
\usepackage[textsize=small]{todonotes}

\ifpdf
  \DeclareGraphicsExtensions{.eps,.pdf,.png,.jpg}
\else
  \DeclareGraphicsExtensions{.eps}
\fi

\newsiamremark{remark}{Remark}
\newsiamremark{hypothesis}{Hypothesis}
\crefname{hypothesis}{Hypothesis}{Hypotheses}
\newsiamthm{claim}{Claim}
\newsiamthm{conjecture}{Conjecture}
\newsiamthm{observation}{Observation}
\newsiamthm{problem}{Problem}
\newsiamthm{question}{Question}
\newsiamthm{assumption}{Assumption}

\crefname{theorem}{theorem}{theorems}
\Crefname{theorem}{Theorem}{Theorems}

\crefname{proposition}{proposition}{propositions}
\Crefname{proposition}{Proposition}{Propositions}

\crefname{lemma}{lemma}{lemmas}
\Crefname{lemma}{Lemma}{Lemmas}

\crefname{corollary}{corollary}{corollaries}
\Crefname{corollary}{Corollary}{Corollaries}

\crefname{definition}{definition}{definitions}
\Crefname{definition}{Definition}{Definitions}

\crefname{claim}{claim}{claims}
\Crefname{claim}{Claim}{Claims}

\crefname{remark}{remark}{remarks}
\Crefname{remark}{Remark}{Remarks}

\crefname{conjecture}{conjecture}{conjectures}
\Crefname{conjecture}{Conjecture}{Conjectures}

\crefname{observation}{observation}{observations}
\Crefname{observation}{Observation}{Observations}

\crefname{problem}{problem}{problems}
\Crefname{problem}{Problem}{Problems}

\crefname{question}{question}{questions}
\Crefname{question}{Question}{Questions}

\crefname{assumption}{assumption}{assumptions}
\Crefname{assumption}{Assumption}{Assumptions}

\newcommand{\identity}{\mathbf{I}}

\newcommand{\Z}{\mathbb{Z}}

\newcommand{\Q}{\mathbb{Q}}
\newcommand{\zero}{\mathbf{0}}
\newcommand{\one}{\mathbf{1}}

\DeclareMathOperator{\slacklattice}{\Lambda}
\DeclareMathOperator{\slackspace}{\Phi}
\newcommand{\mis}{\textsc{Maximum Weighted Independent Set}\xspace}

\makeatletter
\NewCommandCopy\@@pmod\pmod
\DeclareRobustCommand{\pmod}{\@ifstar\@pmods\@@pmod}
\def\@pmods#1{\mkern4mu({\operator@font mod}\mkern 6mu#1)}
\makeatother

\title{\Large Odd faces as a parameter for maximum independent set in embedded graphs}
    \author{Samuel Fiorini\thanks{Department of Mathematics, Université libre de Bruxelles, Brussels, Belgium. Supported in parts by \emph{Fonds de la Recherche Scientifique} - FNRS through research project BD-DELTA-3 (PDR 40028812, 2024--27), and by the \emph{King Baudouin Foundation} through project BD-DELTA-2 (convention 2023-F2150080-233051, 2023--26). (\email{Samuel.Fiorini@ulb.be}).}
    \and Stefan Kober\thanks{Department of Mathematics, London School of Economics and Political Science, London, UK. Part of this work was completed while Stefan Kober was affiliated with Université \mbox{Libre} de Bruxelles, Belgium, and received support by the \emph{Fonds de la Recherche Scientifique} - FNRS through research project BD-DELTA-3 (PDR 40028812).
    (\email{S.Kober@lse.ac.uk}, \url{https://sites.google.com/view/stefankober}).}   
    \and Caleb McFarland\thanks{School of Mathematics, Georgia Institute of Technology, Atlanta, USA. Supported in part by the Georgia Tech ARC-ACO Fellowship and in part by the National Science Foundation under Grant No.~DMS-2452111.
  (\email{cmcfarland30@gatech.edu}, \url{https://sites.google.com/view/caleb-mcfarland/}).}}

\date{}

\begin{document}

\maketitle


\fancyfoot[R]{\scriptsize{Copyright \textcopyright\ 2027 by SIAM\\
Unauthorized reproduction of this article is prohibited}}





\begin{abstract}
    We give a polynomial-time algorithm for finding a maximum weight independent set in a graph embedded in a surface of bounded genus and with a bounded number of vertex-disjoint odd faces. In particular, we obtain an algorithm for graphs with an even-faced embedding in a torus and for planar graphs with a bounded number of odd faces, answering a question of Thomassé (2018) and a question of Choi, Gorsky, Kim, McFarland, and Wiederrecht (2025) respectively. We reduce the original problem to finding a minimum cost non-negative integer circulation with homology constraints in the dual of a certain double cover of the original graph.
\end{abstract}

\section{Introduction}
\subsection{Contribution}

The \mis problem (MIS) takes as input a graph $G$ and a weighting $w:V(G) \rightarrow \mathbb{Q}_{\geq 0}$ and asks for an independent set of maximum weight. MIS is well known to be \textsf{NP}-hard in general, and remains \textsf{NP}-hard even in planar graphs of maximum degree 3 and arbitrarily large girth~\cite{GareyJ1977Rectilinear, Poljak1974ANote}. Due to the problem's hardness in general, there has been a long line of work on determining for which graph classes MIS admits a polynomial time algorithm (see \Cref{subsection:related_work}). A classical example of such a graph class is bipartite graphs, and hence it is natural to ask whether MIS remains tractable on graphs which are ``close'' to being bipartite. One measure of how close a graph is to being bipartite is the \emph{odd cycle packing number} of a graph, which measures the number of pairwise vertex disjoint odd cycles. Recently, Fiorini, Joret, Weltge, and Yuditsky~\cite{fiorini2025integer} showed that MIS admits a polynomial time algorithm on graphs of odd cycle packing number at most $k$ for every constant $k$. Their algorithm is fundamentally based on the work of Conforti, Fiorini, Huynh, Joret, and Weltge~\cite{conforti2020stable} on graphs with bounded genus and odd cycle packing number. Our main result extends the result of Conforti et al.\ to graphs that embed in a surface of bounded genus with a bounded number of vertex disjoint odd \emph{faces}.


\begin{theorem}\label{thm:main_result}
    Let $g,k \in \mathbb{Z}_{\geq 0}$ be fixed. There  exists a polynomial time algorithm that, given a graph $G$ together with a cellular embedding in a surface of Euler genus $g$ with at most $k$ vertex-disjoint odd faces and vertex weights $w: V(G) \rightarrow \Q_{\ge 0}$, finds a maximum weight independent set in $G$ in time $n^{f(g+k)}$ for some computable function $f$.
\end{theorem}

In particular, \Cref{thm:main_result} gives a polynomial time algorithm for MIS in a graph with an even-faced embedding in a torus, asked first by Thomass\'{e}~\cite[Problem~8]{spirkl2018barbados} for the toroidal grid and more generally by Choi, Gorsky, Kim, McFarland, and Wiederrecht~\cite{choi2025odd} for an even-faced embedding.

We also obtain the following as a corollary, answering a question of Choi et al.~\cite{choi2025odd}.
\begin{corollary}
    Let $k \in \mathbb{Z}_{\geq 0}$ be fixed. There exists a polynomial time algorithm that, given a planar graph $G$ with at most $k$ odd faces, finds a maximum weight independent set in $G$.
\end{corollary}

Before our work, the only cases that were known to be polynomially solvable were $k=0$ (bipartite graphs) and $k=2$ (by a result of Gerards~\cite{Gerards1990Graphs}); note the number of odd faces must be even. We provide a more detailed comparison to Gerards' result in \Cref{subsection:related_work}.
We note that planar graphs with two odd faces already have unbounded odd cycle packing number, as do graphs with an even-faced embedding in the torus. Hence our algorithm extends to graphs which are far from bipartite by previous metrics, yet for which the odd cycles have some topological structure. These types of topological constraints on odd cycles arise naturally in the study of odd-minors, and we believe that our results are central to solving MIS in graphs with a forbidden planar odd-minor. We discuss this motivation in more detail below.

We also note that \Cref{thm:main_result} immediately implies the main result of Conforti et al.~\cite{conforti2020stable}, and because our result only relies on an Erd\H{o}s-Pos\'{a} property for odd faces (see \Cref{subsection:removing_odd_faces}) instead of two-sided odd cycles, our analysis is much simpler.


By applying the reduction given in \cite[Section 3]{fiorini2025integer} (see also \cite[Theorem 10.3]{choi2025oddarXiv}), we obtain the following corollary to integer programming. A matrix $A \in \{-1,0,1\}^{m \times n}$ with two nonzero entries per row is the edge-vertex incidence matrix of a signed graph $(G,\gamma)$ if $A$ is supported on the incidence matrix of $G$ and an edge is odd if and only if the two entries in the corresponding row have the same sign. See \Cref{subsection:signed_graphs} for definitions related to signed graphs. We note that the reduction given in \cite[Theorem 10.3]{choi2025oddarXiv} applies to classes of signed graphs closed under subgraphs and parity-preserving subdivisions. The class of graphs with a cellular embedding in a surface of Euler genus $g$ with at most $k$ vertex-disjoint odd faces is closed under parity-preserving subdivisions, but it is not closed under subgraphs because the embedding must be cellular. One can get around this by setting the weights of vertices to zero instead of deleting them.

\begin{corollary}\label{cor:ip}
    Let $g,k \in \mathbb{Z}_{\geq 0}$ be fixed. There exists a polynomial time algorithm to solve the integer program
    $$\max \{w^Tx : Ax \leq b, x \in [\ell,u] \cap \mathbb{Z}^n\}$$
    when $A \in \{-1,0,1\}^{m \times n}$ is the edge-vertex incidence matrix of a signed graph with a given cellular embedding in a surface of Euler genus $g$ with at most $k$ vertex-disjoint odd faces, $\ell,u,w\in\mathbb Z^n$ and $b\in\mathbb Z^m$.
\end{corollary}

Our results are motivated by the study of odd-minors. A graph $H$ is an \emph{odd-minor} of a graph $G$ if $H$ can be obtained from $G$ by deleting vertices, deleting edges, and contracting edge cuts (see \Cref{subsection:signed_graphs} for an equivalent definition). Because every cycle shares an even number of edges with every cut, every cycle in $H$ corresponds naturally to a cycle in $G$ of the same parity. For example, $K_3$-odd-minor-free graphs are exactly bipartite graphs. Odd-minor-free graphs seek to generalize both minor-free graphs and bipartite graphs. They were originally introduced in the context of coloring~\cite{Jensen1995GraphColoring, geelen2009odd, steiner2022asymptotic, kuhn2025disproof}, but have since found applications to algorithmic problems which become tractable on bipartite graphs. To put our results in this context, we note that the class of graphs which embed on a fixed surface with at most $k$ vertex-disjoint odd faces is an odd-minor closed class, see \Cref{lem:even-faced_is_odd_minor_closed}.

For minor-closed classes, an exact algorithmic threshold for MIS is known (we will assume throughout our discussion that $\mathsf{P} \neq \mathsf{NP}$). It is implied by the well-known Grid Theorem of Robertson and Seymour~\cite{RobertsonS1986Graph5} that MIS admits a polynomial-time algorithm on a minor-closed class of graphs if and only if the class forbids a planar graph. Such classes have bounded treewidth, so the structure of such classes is relatively simple and well understood. Choi et al.~\cite{choi2025odd} ask for an analogous algorithmic threshold for odd-minor closed classes, which we conjecture to again be characterized by planar graphs.

\begin{conjecture}\label{conj:planar_odd_minor_tractability}
    Let $\mathcal{C}$ be an odd-minor closed class of graphs. Then \mis admits a polynomial-time algorithm on $\mathcal{C}$ if and only if $\mathcal{C}$ excludes a planar graph.
\end{conjecture}

Unlike graphs with a forbidden planar minor, graphs with forbidden planar odd-minor have a richer and less understood structure, see for instance \cite{choi2025odd}. We expect graphs with an embedding on a surface of bounded genus with a bounded number of odd faces to play an important role in the structure of such graphs, analogous to the role bounded genus graphs play in the structure of $K_t$-minor-free graphs. To see why this should be the case, we note that graphs with an embedding on a surface of bounded genus with a bounded number of odd faces give high-treewidth examples of graphs forbidding a planar odd-minor. Therefore such graphs do not decompose into simpler graphs via standard graph-minors machinery. Previous work also points to the importance of embedded graphs with few odd faces, see for instance \cite{Gerards1990Graphs, fiorini2025integer,choi2025odd}.

\subsection{Related work}\label{subsection:related_work}

We now discuss related work on the \mis (MIS) problem with a particular focus on odd-minor-closed graph classes. Recall that the \emph{odd cycle packing number} (OCP) of a graph $G$ is the maximum number of pairwise vertex-disjoint odd cycles, and the \emph{odd cycle transversal number} (OCT) of a graph is the size of the smallest set of vertices whose deletion leaves the graph bipartite. Clearly a graph has OCP less than $k$ if and only if it forbids $k\cdot K_3$ as an odd-minor. Reed~\cite{Reed1999Mangoes} showed that a graph has bounded OCT if and only if it forbids $k\cdot K_3$ and the Escher wall of order $k$ as an odd-minor. Clearly MIS admits a polynomial time algorithm on graphs of bounded OCT by guessing a hitting set for the odd cycles, but in fact this can be made fixed-parameter tractable by a result of Reed, Smith, and Vetta~\cite{ReedSV2004Finding}.

The problem of MIS on graphs of bounded OCP was motivated by a famous conjecture in integer programming. It is conjectured that any integer program whose coefficient matrix has all subdeterminants bounded by $\Delta$ in absolute value admits a polynomial-time algorithm when $\Delta$ is constant~\cite{shevchenko1996qualitative}. An important special case of this conjecture is the edge formulation of MIS as an integer program, and in this case the conjecture would imply that graphs of bounded OCP admit a polynomial-time algorithm for MIS. The case of OCP equal to 0 is the case of bipartite graphs, while the case of OCP equal to 1 was first solved by Artmann, Weismantel, and Zenklusen~\cite{ArtmannWZ2017Strongly}, who showed more generally that any integer program whose coefficient matrix has all subdeterminants bounded by $2$ in absolute value can be solved in polynomial time.
Conforti, Fiorini, Huynh, Joret, and Weltge~\cite{conforti2020stable} showed that MIS is tractable on graphs with bounded genus and bounded OCP. After applying an Erd\H{o}s-Pos\'{a} type result, the crucial step of their algorithm is solving MIS for graphs with an even-faced embedding on a surface of bounded genus such that the odd cycles are exactly the 1-sided cycles. This is a natural topological structure for odd cycles dating back to Reed's work on the Escher wall~\cite{Reed1999Mangoes}. In our work, we apply a simpler Erd\H{o}s-Pos\'{a} type result to reduce to MIS on graphs with an even-faced embedding on a surface of bounded genus, but now with no further restriction on the odd cycles. Fiorini, Joret, Weltge, and Yuditsky~\cite{fiorini2025integer} extended the work of Conforti et al.~\cite{conforti2020stable} to solve MIS on all graphs of bounded OCP. They did so by showing that graphs of bounded OCP ``almost embed'' on a surface of bounded genus such that odd cycles are exactly the 1-sided cycles.

A different approach to solving MIS in odd-minor-free classes has been the use of width parameters, which largely seek to generalize that MIS admits a polynomial time algorithm on bipartite graphs and graphs of bounded treewidth. Let $\mathcal{B}$ denote the class of bipartite graphs. Examples of such parameters are \emph{elimination distance to $\mathcal{B}$}~\cite{BulianD2016Graph, BulianD2017FixedParameter} and $\mathcal{B}$-treewidth~\cite{EibenGHK2021Measuring, JansendK2021FPTAlgorithms, JansendKW2021Vertex}, which are the minimum $k$ such that a graph $G$ has a set $X \subseteq V(G)$ with $G - X \in \mathcal{B}$ and with the torso of $X$ in $G$ having treedepth (in the case of elimination distance to $\mathcal{B}$) or treewidth (in the case of $\mathcal{B}$-treewidth) at most $k$. Another example is \emph{$\mathcal{B}$-blind-treewidth}~\cite{GollinW2023OddMinors}, which measures the largest treewidth over all non-bipartite blocks of $G$. Jaffke, Morelle, Sau, and Thilikos~\cite{jaffke2025dynamic} introduced \emph{bipartite treewidth}, which gives a proper generalization of both $\mathcal{B}$-blind-treewidth and $\mathcal{B}$-treewidth. A graph has bipartite treewidth at most $k$ if it admits a tree decomposition such that every bag has OCT at most $k$ and the intersection of a bag with its neighbors has size at most 1 after deleting this hitting set. Finally, Choi, Gorsky, Kim, McFarland, and Wiederrecht~\cite{choi2025odd} introduced the notion of \emph{OCP-treewidth}, which generalizes both bipartite treewidth and graphs of bounded OCP. A graph has OCP-treewidth at most $k$ if it admits a tree decomposition such that all adhesions have size at most $k$ and all bags have OCP at most $k$, together with a further restriction that does not allow odd cycles to cross between bags outside a small hitting set. The definition of OCP-treewidth does not require adhesions to have size 1 outside a hitting set as in the definition of bipartite treewidth, but there is a notion of a tame decomposition that has such a constraint, which turns out to be useful both structurally and algorithmically. Choi et al.\ also give a forbidden odd-minor characterization of classes that have bounded OCP-treewidth. They show that an odd-minor closed class has bounded OCP treewidth if and only if it excludes a planar graph with two odd faces and excludes a planar graph with all odd faces sharing a vertex with the outer face.

Separate from the approaches mentioned thus far is the result of Gerards~\cite{gerards1989min}, who proved that MIS admits a polynomial time algorithm on graphs with no $K_4$ odd-minor. This implies that MIS is tractable on planar graphs with 2 odd faces. In fact, planar graphs with 2 odd faces are the key class in his theorem, as Gerards proves in his thesis~\cite{Gerards1990Graphs} that graphs with no $K_4$ odd-minor decompose into parts which are either planar with 2 odd faces, have OCT at most 1, or are equal to a single fixed graph. As previously noted, planar graphs with 2 odd faces can have unbounded OCP, and even unbounded OCP-treewidth. Since 1989 when Gerards' published his result, there have been no generalizations even to planar graphs with 4 odd faces (recall that the number of odd faces must be even). We generalize his result to planar graphs with any constant number of odd faces, and perhaps just as importantly, we give an independent proof of his result using techniques developed by Conforti et al.~\cite{conforti2020stable} which more easily allow for generalizations to odd-minor-free classes as in \cite{fiorini2025integer}. On the surface Gerards' approach is similar to that of Conforti et al.; he finds a special orientation of the graph and uses this orientation to set up a minimum cost circulation. However, Gerards' algorithm is a ``dual'' algorithm in the sense that it computes the maximum weight of an independent set and then uses this to obtain an optimal independent set. In comparison, our algorithm qualifies more as a ``primal'' algorithm; we find the optimal independent set directly.

Our approach gives some hints as to why Gerards' result~\cite{gerards1989min} was not extended to planar graphs with more odd faces before. One way to use our approach to solve MIS in planar graphs with two odd faces is to compute an minimum cost $(s,t)$-path in the dual of the \emph{bipartite double cover} (see \Cref{section:derived_graph}) which in this case will be planar. Hence the minimum cost flow problem will not have any extra topological constraints. However in the case of four or more odd faces, the bipartite double cover will embed on a surface of higher genus, and topological constraints on the flow problem become necessary.

While MIS is $\mathsf{NP}$-hard on any odd-minor closed class containing all planar graphs, Tazari~\cite{Tazari2012Faster} showed that MIS admits a polynomial time approximation scheme on $H$-odd-minor-free graph classes for every graph $H$. This is in stark contrast to the general case, where it is $\mathsf{NP}$-hard to approximate the maximum size of an independent set to within a factor of $n^{1-\varepsilon}$ for every $\varepsilon > 0$~\cite{haastad1996clique, Khot2001Improved}.

The MIS problem has also been studied in graph classes which are closed under induced subgraphs \cite{grotschel1984polynomial, minty1980maximal, lozin2008polynomial, chudnovsky2025treeIV, chudnovsky2025treeV}, induced minors~\cite{Yolov2018Minor, DallardMS2024TreeIndependence, dallard2024treewidth, chudnovsky2024tree, bonnet2023maximum, korhonen2023grid}, and vertex minors~\cite{Oum2005Rankwidth, GeelenKMW2023Grid}. Dallard, Milani\v{c}, and \v{S}torgel asked whether \Cref{conj:planar_odd_minor_tractability} holds when ``odd-minor closed'' is replaced by ``induced minor closed''~\cite{dallard2024treewidth}. This was recently disproven by Bonnet and Chang~\cite{bonnet2026max}. Korhonen~\cite{korhonen2023grid} showed that MIS admits a subexponential time algorithm on the class of $H$-induced-minor-free graphs for every planar graph $H$. For vertex-minors on the other hand, Geelen~\cite{McCarty2021Local} conjectured that MIS admits a polynomial-time algorithm on every proper vertex-minor closed class.

Graphs admitting an even-faced embedding have previously been studied in the context of vertex coloring~\cite{devos2005coloring, kral2008coloring, hutchinson1995three} and face coloring~\cite{mohar2002coloring, hutchinson2002colouring, kral2012coloring, brianski2024four}, in the context of complete odd-minors~\cite{fijavvz2016odd}, and in the context of ribbon graphs~\cite{metsidik2021eulerian, metsidik2020eulerian}.

\subsection{Extensions}
There is still much to be learned about even-faced embeddings of graphs. For instance, can one recognize when a graph has an even-faced embedding in a surface of small genus?

\begin{conjecture}
    For each $g \geq 0$, there exists an algorithm running in time $f(g)n^{\mathcal{O}(1)}$ for some computable function $f$ which takes as input a graph $G$ and either returns an even-faced embedding in a surface of genus $g$ or correctly identifies that no such embedding exists.
\end{conjecture}

If so, this would give an ``even-faced'' analogue of a classic result of Mohar~\cite{mohar1999linear} (see also \cite{kawarabayashi2008simpler}), and would allow us to drop the condition that the embedding of $G$ is given as input in \Cref{thm:main_result}. Even a runtime of the form $n^{f(g)}$ would be interesting. We note that the signed graph minor relation, and hence the odd-minor relation, is expected to be a well-quasi-order; this was announced by Geelen, Gerards, and Whittle~\cite{GeelenGW2014RotasConj} as part of their Matroid Minors Project. It is also expected that one can test for a fixed odd-minor in fixed-parameter tractable time; this should follow from machinery developed by Huynh~\cite{Huynh2009Thesis} and Liu and Yoo~\cite{LiuY2025disjoint}. In particular, Huynh~\cite{Huynh2009Thesis} gives an $\mathsf{XP}$ algorithm in his thesis. Together, this would imply a nonconstructive algorithm for the decision version of the problem.

We conjecture that our algorithm can be made fixed-parameter tractable (FPT). By \Cref{subsection:removing_odd_faces}, it suffices to give an algorithm for even-face embedded graphs.

\begin{conjecture}
    Let $g \in \mathbb{Z}_{\geq 0}$ be fixed. There exists an algorithm that, given a graph $G$ together with a cellular even-face embedding in a surface of Euler genus $g$ and vertex weights $w: V(G) \rightarrow \mathbb{Q}_{\geq 0}$, finds a maximum weight independent set in $G$ in time $f(g)n^{\mathcal{O}(1)}$ for some computable function $f$.
\end{conjecture}

Following this paper's acceptance to SODA 2027, Barish~\cite{Barish2026FPT} claimed a proof of the above conjecture in a Lean repository without an accompanying writeup. Because the case when the odd cycles are exactly the 1-sided cycles reduces to finding a minimum cost flow in a particular homology class, we expect that any answer to the above would also solve the following, which asks for an FPT algorithm to the problem considered in \cite{morell2021minimum}.

\begin{conjecture}
    Let $g \in \mathbb{Z}_{\geq 0}$ be fixed. There exists an algorithm that, given a directed graph $D$ together with a cellular embedding in a surface of Euler genus $g$, edge costs $c: E(G) \rightarrow \mathbb{Q}_{\geq 0}$, and an integer circulation $x$ on $D$, finds a minimum cost integer circulation homologous to $x$ in time $f(g)n^{\mathcal{O}(1)}$ for some computable function $f$.
\end{conjecture}

\section{Algorithm overview}\label{section:algorithm_overview}
We first provide an overview of our algorithm. We will assume familiarity with surfaces and graph embeddings in this section, see \Cref{subsection:embeddings} for relevant background. In particular, we assume that the graph comes together with an embedding scheme $(\pi,\lambda)$. Recall that we consider the following problem.

\begin{problem}\label{problem:MIS_few_odd_faces}
    Given a graph $G$ with a cellular embedding in a surface of Euler genus $g$ with at most $k$ vertex-disjoint odd faces and vertex weights $w: V(G) \rightarrow \mathbb{Q}_{\geq 0}$, compute a maximum $w$-weight independent set in $G$.
\end{problem}

We first reduce to the case when $G$ has an even-faced embedding in a surface of slightly higher genus in \Cref{subsection:removing_odd_faces}. Using results from \cite{Gollin2026}, we show that we may delete at most $\mathcal{O}(k^2g)$ vertices to reduce to the case when $G$ has a constant number of odd faces. Then we may add a handle between odd faces to increase the genus while removing odd faces. We must do these operations while preserving that the graph has a cellular embedding, so instead of deleting vertices we introduce a small gadget with the same effect, see \Cref{subsection:removing_odd_faces} for details. Hence we reduce to the following problem.

\begin{problem}\label{problem:MIS_even_faced}
    Given a graph $G$ with a cellular even-faced embedding in a surface of Euler genus $g$ and vertex weights $w: V(G) \rightarrow \mathbb{Q}_{\geq 0}$, compute a maximum $w$-weight independent set in $G$.
\end{problem}

In \Cref{subsection:standard_assumptions} we apply some preprocessing to reduce to a more structured instance. We summarize the desired properties below. (The notation $\delta(v)$ is defined in \Cref{subsection:definitions} below, together with other preliminaries.)

\begin{definition}
    Let $G$ be a graph. A vertex weighting $w: V(G) \rightarrow \mathbb{Q}_{\geq 0}$ is called \emph{edge-induced} if there is a nonnegative edge weighting $c: E(G) \rightarrow \mathbb{Q}_{\geq 0}$ such that $w(v) = \sum_{e \in \delta(v)} c(e)$ for every vertex $v \in V(G)$.
\end{definition}

\begin{definition}
    A weighted graph $(G,w)$ is said to satisfy the \emph{standard assumptions} if $G$ is 2-connected, non-bipartite, and even-faced embedded in a surface of Euler genus $g$ with facewidth at least 2, and the weighting $w$ is edge-induced.
\end{definition}

Because our weights are edge-induced by $c: E(G) \rightarrow \mathbb{Q}_{\geq 0}$, which we think of as the \emph{cost} of an edge, it now suffices to find a minimum cost set of edges which are given by the uncovered edges of some independent set. We call such edges \emph{slack}. We note already that edge-induced weights intrinsically encode that the problem is easy for bipartite graphs, as for bipartite graphs choosing one side of the bipartition will have maximum total weight.

In \Cref{subsection:slack_vectors} we recall the results of Conforti et al.~\cite{conforti2020stable} to precisely describe which sets of edges appear as slack edges in an independent set. That is, if $A \in \{0,1\}^{E(G) \times V(G)}$ is the edge-vertex incidence matrix of $G$, then $y \in \mathbb{Z}_{\geq 0}^{E(G)}$ is called a \emph{slack vector} if there exists $x \in \mathbb{Z}^{V(G)}$ such that $Ax + y = \mathbf{1}$. For $C = (v_0,e_1,v_1,e_2,v_2,\dots, v_{k-1},e_k,v_0)$ a closed walk in $G$ and $y \in \mathbb{Z}^{E(G)}$, we define
$$\omega_C(y) := \sum_{i=1}^k (-1)^{i-1}y(e_i).$$
Conforti et al.\ showed that $y \in \mathbb{Z}_{\geq 0}^{E(G)}$ is a slack vector if and only if $\omega_C(y) = 0$ for all even closed walks and $\omega_C(y) \equiv 1 \pmod 2$ for all odd closed walks. We show that when $G$ is 2-connected and non-bipartite, there is a collection of $|E(G)| - |V(G)| + 1$ cycles whose constraints imply the above constraints for every closed walk, and hence we reduce to the below problem.

\begin{restatable}{problem}{MinCostFeasibleSlack}\label{problem:min_cost_feasible_slack}
    Given a graph $G$ satisfying the standard assumptions, non-negative edge costs $c: E(G) \rightarrow \mathbb{Q}_{\geq 0}$, and a collection $\mathcal{C}$ of $|E(G)| - |V(G)| + 1$ cycles with exactly one of them odd, minimize $c(y) = \sum_{e \in E(G)} c(e)y(e)$ over all $y \in \mathbb{Z}_{\geq 0}^{E(G)}$ satisfying $\omega_C(y) = 0$ for every even cycle $C \in \mathcal{C}$ and $\omega_C(y) \equiv 1 \pmod 2$ for the unique odd cycle $C \in \mathcal{C}$. Furthermore, $\mathcal{C}$ is such that the above constraints exactly characterize the slack vectors $y \in \mathbb{Z}^{E(G)}_{\geq 0}$.
\end{restatable}

In the work of Conforti et al.~\cite{conforti2020stable}, the authors then find a special orientation of the dual graph such that the solution $y \in \mathbb{Z}_{\geq 0}^{E(G)}$ can be interpreted as a circulation in this orientation. This has the advantage that the constraint $\omega_C(y) = 0$ for facial cycles of $G$ is implied by flow conservation in this dual graph. The type of orientation needed is an \emph{alternating orientation}. An orientation of an embedded graph $G^*$ is called \emph{alternating} if for each vertex $v \in V(G^*)$, edges alternatively leave and enter $v$ in the local cyclic order of edges incident to $v$. Finding such an alternating orientation of $G^*$ is only possible in the setting of \cite{conforti2020stable} (see \Cref{lem:alternating_orientation_equivalence}).

To get around this problem, we consider the \emph{derived graph} $\tilde{G}$ of $G$ (for a particular signing of $G$), see \Cref{section:derived_graph}. The derived graph $\tilde{G}$ is constructed as follows. The graph $G$ is equipped with an embedding scheme $(\pi, \lambda)$, where $\pi$ gives the local cyclic orientation around each vertex and $\lambda: E(G) \rightarrow \mathbb{Z}_2$ encodes whether an edge is ``twisted'', see \Cref{subsection:embeddings}. The vertex set of $\tilde{G}$ consist of two copies of $V(G)$, and we denote for each vertex $v \in V(G)$ its two copies by $v_0,v_1$. For each edge $uv \in E(G)$ with $\lambda(uv) = 0$, one adds two edges in $E(\tilde{G})$ between opposite copies of $u$ and $v$, that is, one adds the edges $u_0v_1$ and $u_1v_0$. For each edge $uv \in E(G)$ with $\lambda(uv) = 1$, one adds two edges in $E(\tilde{G})$ between the same copies of $u$ and $v$, that is, one adds the edges $u_0v_0$ and $u_1v_1$. It turns out that $\tilde{G}$ embeds in a surface of genus $2g - 2$, and that its dual $\tilde{G}^*$ has an alternating orientation. We defer the reader to \Cref{section:derived_graph} for a topological construction of $\tilde{G}$.

We now have an alternating orientation $D^*$ of $\tilde{G}^*$. We map circulations $f \in \mathbb{Z}_{\geq 0}^{E(D^*)}$ to vectors $y \in \mathbb{Z}_{\geq 0}^{E(G)}$ via the natural projection $\rho: \mathbb{Z}^{E(D^*)} \rightarrow \mathbb{Z}^{E(G)}$ where $y = \rho(f)$ is given by $y(e) = f(\vec{e_0}) + f(\vec{e_1})$ for $e_0,e_1 \in E(\tilde{G})$ the copies of $e \in E(G)$. Now flow conservation at each vertex in $D^*$ implies $\omega_C(y) = 0$ for every facial cycle $C$. By Euler's formula, it suffices to add at most $g$ extra constraints $\omega_{C_1}(\rho(f)) \equiv 1 \pmod 2$ for odd cycle $C_1$ and $\omega_{C_i}(\rho(f)) = 0$ for even cycles $C_2, \dots, C_g$. Somewhat surprisingly, we show that an optimal solution $y \in \mathbb{Z}^{E(G)}$ to \Cref{problem:min_cost_feasible_slack} can be found as the image of a minimum cost circulation $f$ in $D^*$ satisfying the above constraints, and hence we reduce to the following problem in \Cref{thm:reduction_to_integer_circulation}.

\begin{restatable}{problem}{IntegerCirculation}\label{problem:integer_circulation}
    Let $D^*$ be an alternating orientation of $\tilde{G}^*$ for $G$ satisfying the standard assumptions, and let $C_1, \dots, C_g$ be cycles in $G$ where $C_1$ is odd and $C_2, \dots, C_g$ are even. Find a minimum $c$-cost nonnegative integer circulation $f$ in $D^*$ satisfying the linear constraints
    \begin{enumerate}
        \item $\omega_{C_i}(\rho(f)) = 0$ for every even cycle $C_2, \dots, C_g$, and
        \item $\omega_{C_1}(\rho(f)) \equiv 1 \pmod 2$.
    \end{enumerate}
\end{restatable}

Finally, in \Cref{section:final_algorithm}, we show that the constraints on the circulation $f$ in $D^*$ only depend on the $\mathbb{Z}$-homology of $f$ (see \Cref{subsection:homology} for related definitions). This implies that we can solve \Cref{problem:integer_circulation} in polynomial time by using the work of Chambers, Erickson, and Nayyeri~\cite{chambers2012homology} for the orientable case and Morell, Seidel, and Weltge~\cite{morell2021minimum} for the nonorientable case.

\section{Preliminaries}

\subsection{Definitions}\label{subsection:definitions}
All graphs will be simple unless otherwise stated. We use standard graph theory notation. In particular, we use $V(G)$ and $E(G)$ to denote the vertex and edge set of a graph respectively, and we use $\delta(v)$ to denote the set of edges incident to $v \in V(G)$. We will often identify a function $x: A \rightarrow B$ with a vector $x \in B^A$, and for instance use the notation $x(a)$ to denote the $a$-th entry of $x \in B^A$. We use the notation $[k] := \{1, 2, \dots, k\}$.

In \Cref{subsection:removing_odd_faces} we will need the following definition.

\begin{definition}
    A \emph{tree-decomposition} of a graph $G$ is a pair $(T,\beta)$ where $T$ is a tree and $\beta: V(T) \rightarrow 2^{V(G)}$ is such that
    \begin{enumerate}
        \item for every edge $e \in E(G)$, there exists $t \in V(T)$ such that both ends of $e$ are contained in $\beta(t)$, and
        \item for every $v \in V(G)$, the set $\beta^{-1}(v)$ induces a connected, non-empty subtree of $T$.
    \end{enumerate}
    The \emph{width} of a tree decomposition is $\max_{t \in V(T)} |\beta(t)| - 1$. The \emph{treewidth} of a graph $G$ is the minimum width over all tree decompositions of $G$. We call $\beta(t)$ the \emph{bag} of $t$.
\end{definition}

\subsection{Signed graphs}\label{subsection:signed_graphs}
To understand our algorithm it will be most convenient to work with signed graphs. A \emph{signed graph} is a pair $(G,\gamma)$ where $G$ is a graph and $\gamma:E(G) \rightarrow \mathbb{Z}_2$ assigns a parity to every edge. The parity of a cycle $C = v_1v_2\dots v_kv _1$ in $(G,\gamma)$ is given by $\sum_{i =1}^k \gamma(v_iv_{i+1})$ where indices are taken modulo $k$, and the sum is computed in $\Z_2$. In the context of this paper we typically associate an unsigned graph $G$ with the signed graph $(G,\mathbbm{1})$ where $\mathbbm{1}:E(G) \rightarrow \mathbb{Z}_2$ maps every edge to 1.

A signed graph $(G,\gamma')$ is obtained from $(G,\gamma)$ by \emph{shifting at $v \in V(G)$} if $\gamma'(e) = \gamma(e)$ for all $e \in E(G) \setminus \delta(v)$ and $\gamma'(e) = \gamma(e) + 1$ for all $e \in \delta(v)$. More generally $(G,\gamma')$ is obtained from $(G,\gamma)$ by \emph{shifting} if $(G,\gamma')$ is obtained from $(G,\gamma)$ by a sequence of shiftings at vertices. It is well known that two signed graphs have the same collection of odd cycles if and only if one is obtained from the other by shifting.

A signed graph $(G,\gamma_G)$ contains another signed graph $(H,\gamma_H)$ as a \emph{minor} if $(H,\gamma_H)$ can be obtained from $(G,\gamma)$ by shifting, deleting vertices, deleting edges, and contracting even edges. We note that we disallow the contraction of odd-loops; this is sometimes called a \emph{link minor} in the associated matroid. We can use this to give an equivalent definition of odd-minors: $G$ contains $H$ as an \emph{odd-minor} if and only if the signed graph $(G,\mathbbm{1})$ contains the signed graph $(H,\mathbbm{1})$ as a minor.

In the below, the parity of a face in a signed graph $(G,\gamma)$ is the sum of $\gamma(e)$ over all edges $e$ on its boundary, where edges that appear twice along the boundary are counted twice. The parity of a face in an unsigned graph $G$ is the parity of the face in $(G,\mathbbm{1})$.
\begin{lemma}\label{lem:even-faced_is_odd_minor_closed}
    Let $\mathcal{C}$ denote the class of graphs that admit an embedding in a surface $\Sigma$ with at most $k$ vertex-disjoint odd faces. Then $\mathcal{C}$ is odd-minor-closed.
\end{lemma}
\begin{proof}
    This follows from the more general fact that the class $\mathcal{C}'$ of signed graphs with an embedding in $\Sigma$ with at most $k$ vertex-disjoint odd faces is minor-closed. Indeed, $G \in \mathcal{C}$ if and only if $(G,\mathbbm{1}) \in \mathcal{C}'$. Then for any odd-minor $H$ of $G$, $(H,\mathbbm{1})$ is a minor of $(G,\mathbbm{1})$, and hence $(H,\mathbbm{1})$ is in $\mathcal{C}'$.
    
    We then simply note the following properties. Clearly shifting, deleting edges, deleting vertices, and contracting even edges leaves the signed graph embeddable in $\Sigma$, so we argue that the number of vertex-disjoint odd faces does not increase. Shifting at a vertex does not change the parity of any face. Deleting an edge either doesn't change the parity of any face (if that edge appears twice on the boundary of some face), or it merges two faces such that the new face has parity equal to the sum of the previous faces. It suffices to only ever delete isolated vertices, and clearly this does not change the parity of a face. Finally, contracting an even edge is either the same as deleting an edge (if the edge is a loop) or does not affect the parity of any incident faces.
\end{proof}

We note some subtlety here. The class of graphs with a \emph{cellular} embedding in a surface with a bounded number of (vertex-disjoint) odd faces is not odd-minor closed. Deleting vertices and edges may cause the embedding to be non-cellular.

\subsection{Surfaces and graph embeddings}\label{subsection:embeddings}

In this subsection we review the basics of surfaces and combinatorial graph embeddings. We refer the reader to Mohar and Thomassen's book~\cite{MoharT2001Graphs} for further background.

Throughout the paper a \emph{surface} $\Sigma$ is a compact 2-dimensional manifold without boundary. We let $\Sigma(h,c)$ be the surface obtained by cutting $2h+c$ disjoint open disk from the sphere and gluing $h$ handles and $c$ crosscaps to the boundaries of these disks. By the classification theorem for surfaces, every surface is homeomorphic to $\Sigma(h,c)$ for some $h$ and $c$. If $c = 0$ the surface is said to be \emph{orientable}, otherwise the surface is said to be \emph{non-orientable}. The \emph{Euler genus} of a surface $\Sigma \cong \Sigma(h,c)$ is $2h+c$.

A simple closed curve in $\Sigma$ is said to be \emph{1-sided} if it has a neighborhood that is homeomorphic to a M\"{o}bius strip, and \emph{2-sided} if it has a neighborhood that is homeomorphic to a cylinder. A simple, closed curve $C$ in $\Sigma$ is said to be \emph{contractible} if some component of $\Sigma \setminus C$ is an open disk whose boundary is $C$. Otherwise we say $C$ is \emph{non-contractible}.

An \emph{embedding} $\psi$ of a graph $G$ in $\Sigma$ is a map from vertices of $G$ to distinct points of $\Sigma$ and edges of $G$ to internally disjoint simple open curves in $\Sigma$ between the image of its endpoints. We use $\psi(G)$ to denote the union of all points and curves in the image of $\psi$. The \emph{faces} of an embedded graph $G$ are the connected components of $\Sigma \setminus \psi(G)$, and we denote the set of faces by $F(G)$. If every face is homeomorphic to an open disk, then we say $\psi$ is a \emph{cellular embedding}.

A curve in $\Sigma$ is \emph{G-normal} if it only meets $\psi(G)$ in vertices. If $\Sigma$ is not the sphere, the \emph{facewidth} of the embedded graph $G$ is the minimum number of vertices intersected by a $G$-normal, non-contractible, simple closed curve in $\Sigma$. Facewidth is not defined for graphs embedded in the sphere. If $G$ is 2-connected and cellularly embedded in $\Sigma$ with facewidth at least 2, then the boundary of each face is a cycle in $G$~\cite[Proposition 5.5.11]{MoharT2001Graphs}.

If $G$ is connected and cellularly embedded in a surface $\Sigma$ of Euler genus $g$, then
$$|V(G)| - |E(G)| + |F(G)| = 2 - g.$$

Given a cycle $C$ in a graph $G$ with embedding $\psi$ in $\Sigma$, we say that $C$ is \emph{1-sided} (resp. \emph{2-sided}) if the simple closed curve $\psi(C)$ is 1-sided (resp. 2-sided).

A cellular embedding of a graph $G$ on a surface $\Sigma$ is determined by its \emph{embedding scheme} $(\pi,\lambda)$ where $\pi$ associates to each vertex of $G$ a cyclic ordering of its incident edges, which one thinks of as the ``clockwise ordering'' of edges around a vertex, and $\lambda: E(G) \rightarrow \mathbb{Z}_2$ takes value 1 if and only if the orientation given by $\pi$ of the edge's endpoints disagree. Because a cyclic ordering of the edges incident to a vertex of degree at most 2 does not encode a clockwise direction, the above definition only makes sense for graphs of minimum degree 3. To extend the definition to all graphs $G$, we first suppress degree 2 vertices to get a multigraph $G'$ and we define the embedding scheme $(\pi',\lambda')$ as above, where every edge incident to a vertex of degree 1 may have an arbitrary $\lambda'$ label. Now for each edge $e'$ of $G'$ obtained by suppressing degree 2 vertices along a path $P$ in $G$, we assign the value of $\lambda(e)$ for edges $e \in E(P)$ such that the total $\lambda$-weight of $P$ is equal to $\lambda'(e')$. This defines an embedding scheme for $G$. We note a cycle in $G$ is one-sided if and only if it is odd with respect to $\lambda$. All embeddings in this paper will be stored as embedding schemes for algorithmic purposes.

Note that if one reverses the orientation of $\pi$ at a vertex $v \in V(G)$ and adds 1 (over $\Z_2$) to the label $\lambda(e)$ for all non-loop edges $e \in \delta(v)$, then this defines the same embedding. We say an embedding scheme $(\pi',\lambda')$ is obtained from $(\pi,\lambda)$ by \emph{shifting at $v$} if it is obtained in this way. We say that two embeddings $(\pi,\lambda)$ and $(\pi',\lambda')$ of a graph $G$ are \emph{(shifting) equivalent} if $(\pi',\lambda')$ is obtained from $(\pi,\lambda)$ by a sequence of shiftings at vertices. Note that if $(\pi,\lambda)$ and $(\pi',\lambda')$ are shifting equivalent, then the signed graph $(G,\lambda')$ is obtained from $(G,\lambda)$ by shifting.

Given an embedding scheme $(\pi, \lambda)$ of a graph $G$, we can naturally define the set of facial walks as walks which, at every vertex $v$, choose the next edge in the cyclic ordering $\pi(v)$ with a direction that swaps whenever the traversed edge has $\lambda(e) = 1$.

We are now equipped to say exactly when the dual of a graph has an alternating orientation.

\begin{lemma}\label{lem:alternating_orientation_equivalence}
    Let $G$ be cellularly embedded in a surface such that the boundary of every face is a cycle in $G$. Then the following are equivalent.
    \begin{enumerate}
        \item $G^*$ admits an alternating orientation,
        \item A cycle in $G$ is odd if and only if it is one-sided,
        \item The constraints $\{\omega_C(y) = 0 \mid C\text{ is a facial cycle}\}$ are not linearly independent.
    \end{enumerate}
\end{lemma}
\begin{proof}
    Consider an embedding scheme $(\pi, \lambda)$ for $G$. Such an embedding satisfies that the odd cycles are exactly the one-sided cycles if and only if $(\pi, \lambda)$ is shifting equivalent to $(\pi', \mathbbm{1})$.
    That is, we have a local orientation $\pi'$ around each vertex such that for every edge, the orientations of its ends disagree. We may then orient each facial cycle of $G^*$ to be a directed cycle agreeing with the orientation of $\pi'$; such an orientation is an alternating orientation of $G^*$. See \Cref{fig:orientations} (see also \cite[Lemma 17]{conforti2020stable}).
    \begin{figure}[ht]
        \centering
        \includegraphics[width=.4\textwidth]{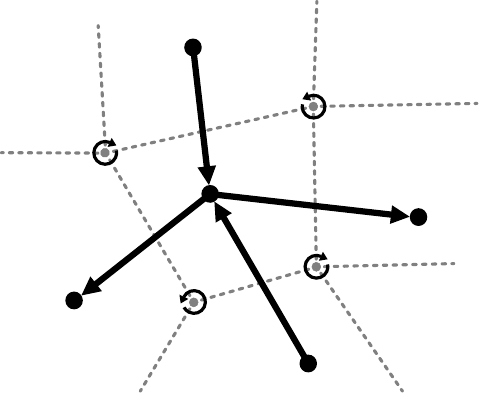}
        \caption{The graph $G$ is shown with dashed edges, and circular arrows displaying the local orientation $\pi'$. Note that neighboring local orientations disagree, which gives an alternating orientation of the dual graph $G^*$ depicted with the solid black arcs.}
        \label{fig:orientations}
    \end{figure}
    
    Conversely, suppose $G^*$ admits an alternating orientation. Then each facial cycle of $G^*$ is a directed cycle in this orientation. This naturally gives a local orientation around the corresponding vertex in $G$, and these orientations have the property that for any edge of $G$, the orientation of its ends disagree. Thus the first two items are equivalent.

    Now suppose $G^*$ admits an alternating orientation. Then, for each facial cycle $C$ corresponding to $v \in V(G^*)$, sign its edges such that edges leaving $v$ get a positive sign and edges entering $v$ get a negative sign. We can possibly negate $\omega_C(y)$ to match this signing, and summing this up over every facial cycle gives a linear dependence between $\omega_C(y) = 0$ for facial cycles $C$.

    Conversely suppose that there exists some linear dependence among $\{\omega_C(y) = 0\}$. Write the constraints as $My = 0$, where $M \in \{-1,0,1\}^{F(G) \times E(G)}$. Note the support of $M$ is the incidence matrix of $G^*$. By assumption, there is some nonzero vector $\mu$ such that $\mu^{\intercal} M = 0$. Since $G^*$ is connected and $\mu \neq 0$, it must be that $\mu$ is nonzero in every entry. Since $M$ has two nonzero entries per column, we may assume $\mu \in \{-1,1\}^{F(G)}$ by scaling $\mu$. By multiplying each row $f$ of $M$ by $\mu(f)$, we may assume $\mu = \mathbbm{1}$. The new matrix $M$ has two nonzero entries per column, and every column contains a 1 and $-1$. Then $M$ defines an alternating orientation of $G^*$ (recall that the rows are scalings of $\omega_C$).
    
\end{proof}

\subsection{Homology and even-faced embeddings}\label{subsection:homology}
Let $(G,\gamma)$ be an embedded signed graph. Given a face with facial walk $(v_0,e_1,v_1,e_2,v_2,\dots, v_{k-1},e_k,v_0)$, the face's \emph{parity} is given by $\sum_{i=1}^k \gamma(e_i) \in \mathbb{Z}_2$. (We define the parity of a face in an unsigned embedded graph by identifying the unsigned graph $G$ with the signed graph $(G,\mathbbm{1})$.) We say a face is \emph{even} or \emph{odd} if its parity is 0 or 1 in $\mathbb{Z}_2$ respectively. We say the embedding is \emph{even-faced} if no odd cycle bounds a disk in $\Sigma$. If the embedding is cellular, this is equivalent to requiring that every face is even.

Let $D$ be a directed graph cellularly embedded in a surface $\Sigma$. Two circulations $f,g \in \mathbb{Z}^{E(D)}$ are said to be \emph{$\mathbb{Z}$-homologous} if $f - g$ is a linear combination of facial circulations with integer coefficients. The \emph{first homology group of $\Sigma$ over the integers}, denoted $H_1(\Sigma; \mathbb{Z})$, is the quotient group of all integer circulations in $D$ by the subgroup of all circulations which are $\mathbb{Z}$-homologous to the zero circulation. We note that the resulting group does not depend on our choice of $D$. It is well known that if $\Sigma$ is orientable of Euler genus $g$ then $H_1(\Sigma; \mathbb{Z}) \cong \mathbb{Z}^g$, and if $\Sigma$ is non-orientable of Euler genus $g$ then $H_1(\Sigma; \mathbb{Z}) \cong \mathbb{Z}_2 \times \mathbb{Z}^{g-1}$.

When $(G,\gamma)$ is a signed graph with an even-faced cellular embedding in a surface $\Sigma$, the parity of a cycle only depends on the ``$\Z$-homology'' of the cycle. We now make this precise. Let $D$ be an arbitrary orientation of $G$. We define the parity of an integer circulation $f$ in $D$ to be $\sum_{e \in E(D)} f(e)\gamma(e) \in \Z_2$. Because $(G,\gamma)$ is even-faced, the parity of any facial circulation is 0. Hence the parity of $f$ depends only on the $\Z$-homology of $f$ (in fact, this defines a homomorphism $\varphi:H_1(\Sigma; \Z) \rightarrow \Z_2$). Then, for any cycle $C$ in $G$, the parity of $C$ is equal to the parity of a $\{-1,0,1\}$-circulation $f$ supported on $E(C)$. Hence the parity of $C$ depends only on the $\Z$-homology of the corresponding circulation. In cases such as this when the direction of flow around $C$ does not matter, we may speak of the ``$\Z$-homology'' of $C$ to mean the $\Z$-homology of the corresponding circulation $f$ in an arbitrary orientation $D$ of $G$.

\section{Preprocessing}

\subsection{Removing odd faces}\label{subsection:removing_odd_faces}
In this section we will reduce from \Cref{problem:MIS_few_odd_faces} to \Cref{problem:MIS_even_faced}. That is, we will reduce to the case when all faces are even by possibly guessing on a small number of vertices and increasing the genus.

We first reduce to the case with a constant number of odd faces, not just a constant number of vertex-disjoint odd faces. For this we use an Erd\H{o}s-Pos\'{a} property for odd faces, which is implied for instance by a more general result due to the authors of~\cite{Gollin2026} which says that Erd\H{o}s-Pos\'{a} holds even when replacing ``odd faces'' by any property for faces. To derive the below statement from \cite[Corollary 4.7]{Gollin2026}, consider coloring every odd face by a unique color and coloring every even face by the same color. We give the proof here for completeness. Recall that a face is odd if its facial walk has odd length.

\begin{theorem}\cite{Gollin2026}\label{thm:odd_faces_EP}
    Let $g,k \in \mathbb{Z}_{\geq 0}$. Let $G$ be a graph embedded on a surface $\Sigma$ of Euler genus $g$. Then either there exist $k$ odd faces $F_1, \dots, F_k$ of $G$ whose boundaries are pairwise vertex-disjoint, or there exists a set $X \subseteq V(G)$ of size at most $\mathcal{O}(k^2g)$ such that every odd face contains some vertex in $X$ on its boundary.

    Moreover, there exists an algorithm running in time $2^{\mathcal{O}(g^2)}\mathsf{poly}(k)|V(G)|^4$ that finds either the faces $F_1, \ldots, F_k$ or the set $X$ with the properties promised above.
\end{theorem}

We will skip the proof of the FPT algorithm which finds $X$. It will become clear from the proof that this reduces to finding $X$ in the low treewidth case, which one can do via dynamic programming. We first show that the above holds for graphs of small treewidth.

\begin{lemma}\cite{Gollin2026}\label{lem:low_tw_odd_faces_EP}
    Let $k \in \mathbb{Z}_{\geq 0}$. Let $G$ be a graph with treewidth at most $\eta$ embedded in a surface $\Sigma$, and let $C_1, \dots, C_m$ be the boundaries of some collection of faces. Then either there exist $k$ vertex disjoint facial boundaries $C_{i_1}, \dots, C_{i_k}$, or there exists a set $X \subseteq V(G)$ of size at most $(k-1)(\eta+1)$ hitting all $C_i$. 
\end{lemma}
\begin{proof}
    Let $(T,\beta)$ be a tree decomposition of $G$ of width at most $\eta$. For each $i \in [m]$, let $T_i$ be the subtree of $T$ induced by all bags containing a vertex of $C_i$. That is, $V(T_i) = \{t \in V(T) : \beta(t) \cap V(C_i) \neq \varnothing\}$. We now show that either there exists $k$ vertex disjoint subtrees $T_{i_1}, \dots, T_{i_k}$ from our collection, or there exists a set $Y \subseteq V(T)$ of size at most $k-1$ hitting all subtrees in our collection. We proceed by induction on $|V(T)|$. Clearly we may assume $|V(T)| \geq 2$, and let $\ell \in V(T)$ be a leaf. If there exists a subtree $T_i$ in our collection consisting only of the vertex $\ell$, then we apply induction to $T - \ell$ with the collection obtained from deleting all subtrees containing $\ell$. If we obtain a packing of disjoint subtrees then we may add $T_i$ to the packing, and if we obtain a hitting set then we may add $\ell$ to the hitting set. If there is no such tree with vertex set equal to $\ell$, then we delete $\ell$ from $T$ and from every subtree of our collection, and again the desired statement then follows by induction.
    
    If we obtain the packing of disjoint subtrees $T_{i_1}, \dots, T_{i_k}$, then the corresponding collection of facial boundaries $C_{i_1}, \dots, C_{i_k}$ are vertex disjoint. If we instead obtain the hitting set $Y \subseteq V(T)$ then the set $X = \{v \in \beta(t) : t \in Y\}$ is a hitting set of size at most $(k-1)(\eta + 1)$ hitting all facial boundaries $C_i$.
\end{proof}

We can now give the proof of \Cref{thm:odd_faces_EP}. This essentially follows from greedily packing odd faces and using the fact that in a graph of bounded genus, any collection of faces that all intersect a common face in at least one vertex induces a graph of bounded treewidth.

\begin{proof}[Proof of \Cref{thm:odd_faces_EP}]
    We first greedily pack odd faces $F_1, \dots, F_b$ whose boundaries are vertex disjoint. We may assume $b < k$, otherwise we obtain the desired collection of faces. Now every other odd face shares some vertex with one of $F_1, \dots, F_b$.

    For each $i \in [b]$, let $G_i$ be the graph induced on the boundary of $F_i$ and all faces which share a vertex with $F_i$. We first claim that $G_i$ has treewidth $\mathcal{O}(g)$. Let $G_i'$ be formed from $G_i$ by adding a vertex $v_F$ inside every face $F$ of $G_i$ which is also a face in $G$ and making $v_F$ incident to every vertex on the facial boundary of $F$. We then note that $G_i'$ has radius at most 3; by definition every vertex is distance at most 3 from $v_{F_i}$. Eppstein \cite[Theorem 2]{Eppstein2000Diameter} proved that graphs of radius $r$ embedded in a surface of Euler genus $g$ have treewidth $\mathcal{O}(rg)$. Because $G_i$ is a subgraph of $G_i'$, we immediately obtain that $G_i$ has treewidth $\mathcal{O}(g)$.

    Now in each $G_i$, we apply \Cref{lem:low_tw_odd_faces_EP} to the collection of odd faces which are also faces of $G$. We may assume each $G_i$ does not contain a packing of $k$ such odd faces, else we obtain a packing in $G$. Thus we obtain a set of $\mathcal{O}(kg)$ vertices $X_i$ hitting all such faces in $G_i$. Then $X = \bigcup_{i = 1}^bX_i$ is the desired hitting set in $G$.
\end{proof}

Now suppose the graph is embedded on a surface of Euler genus $g$ with $k$ vertex-disjoint odd faces. By the above, there exists a set $X$ of at most $\mathcal{O}(k^2g)$ vertices hitting all odd faces. We guess the intersection of a maximum independent set with $X$. Call this intersection $Y \subseteq X$. Guessing $Y$ increases the running time by a factor of $2^{\mathcal{O}(k^2g)}$ which we assume to be constant. One would like to delete $X \cup N_G(Y)$ and find a maximum independent set in the resulting graph, but doing so may destroy that the embedding is cellular. Instead, we do the following. We first set the weight of $X$ and all neighbors of $Y$ to zero. We may assume that no zero-weight vertex is chosen in a maximum weight independent set, and so any optimal solution can be extended to an independent set by including $Y$.

Now, for each $x \in X$, $x$ may be contained in the boundary of many odd faces. However, because $x$ has weight 0, we may subdivide edges incident to $x$ with weight 0 vertices without changing an optimal solution. We now show that, for each $x \in X$, we may subdivide incident edges in such a way that at most one of the faces containing $x$ on the boundary is odd after subdivision. By doing this one at a time for each vertex in $X$, we reduce to an instance with at most $\mathcal{O}(k^2g)$ odd faces which is cellularly embedded.

\begin{figure}[h!]
\centering
\includegraphics[width=.7\textwidth]{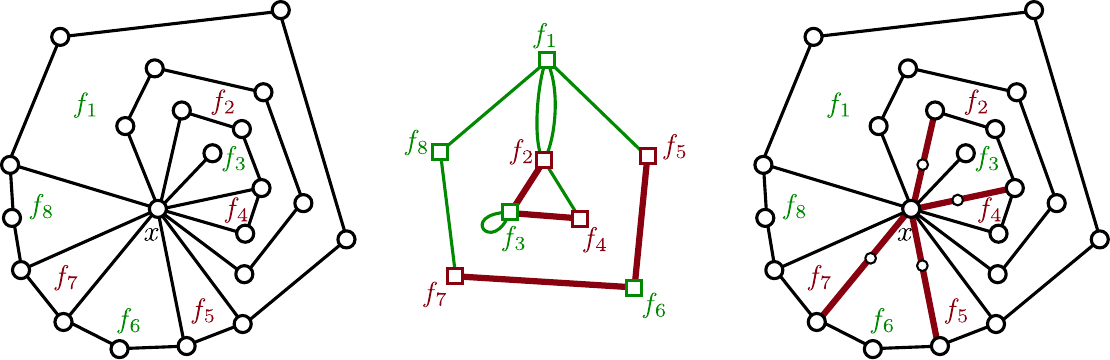}
\caption{The left and right images show the local structure in $G$ around vertex $x$. The image in the center shows the 
subgraph $H_x \subseteq G^*$. In this example, one has $T = T' = \{f_2,f_4,f_5,f_7\}$. The $T'$-join is shown in bold red.}
\label{fig:H_x}
\end{figure}

For each $x \in X$, we construct an auxiliary multigraph $H_x$ as follows, see Figure~\ref{fig:H_x} for an illustration. $H_x$ is a subgraph of the dual $G^*$ consisting of all vertices and edges lying on the face $x$ of $G^*$. That is,
the vertex set of $H_x$ is the set of faces of $G$ with $x$ on their boundary. The edge set of $H_x$ is in bijection with $\delta(x)$, where each edge $e \in \delta(x)$ corresponds to an edge in $H_x$ between the faces containing $e$. We note that $H_x$ is not necessarily a cycle as $x$ may appear multiple times on a facial walk in $G$, and in fact $H_x$ may contain parallel edges or loops. However, we note that $H_x$ is connected. We let $T \subseteq V(H_x)$ be the set of odd faces containing $x$ on their boundary. If $|T|$ is even, we let $T' = T$, otherwise we let $T'$ be formed from $T$ by removing a single element. Because $H_x$ is connected, there exists a set $F \subseteq E(H_x)$ which is a $T'$-join. That is, every vertex $v \in V(H_x)$ has odd $F$-degree if $v \in T'$ and even $F$-degree if $v \in V(H_x) \setminus T'$. Then by subdividing the edges of $\delta(x)$ which correspond to an edge in $F$, we have that every face except possibly the single face in $T \setminus T'$ is even in the resulting graph.

\begin{figure}[ht]
\centering
\includegraphics[width=.45\textwidth]{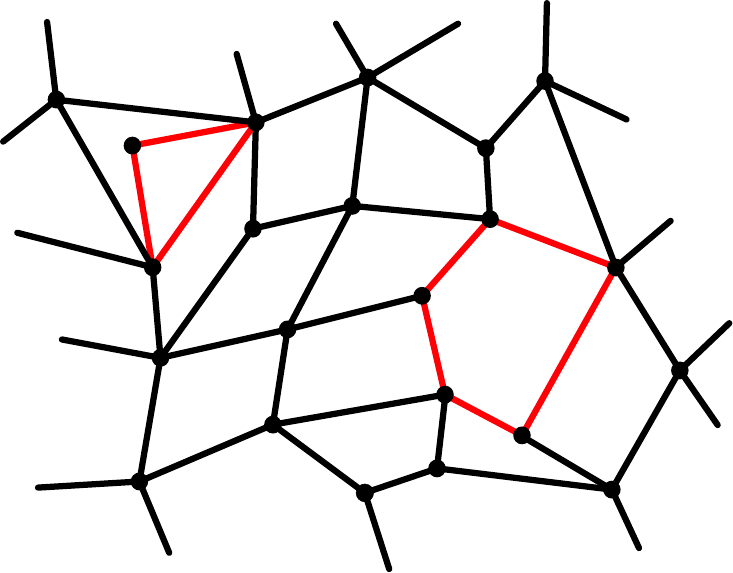}\hfill
\includegraphics[width=.45\textwidth]{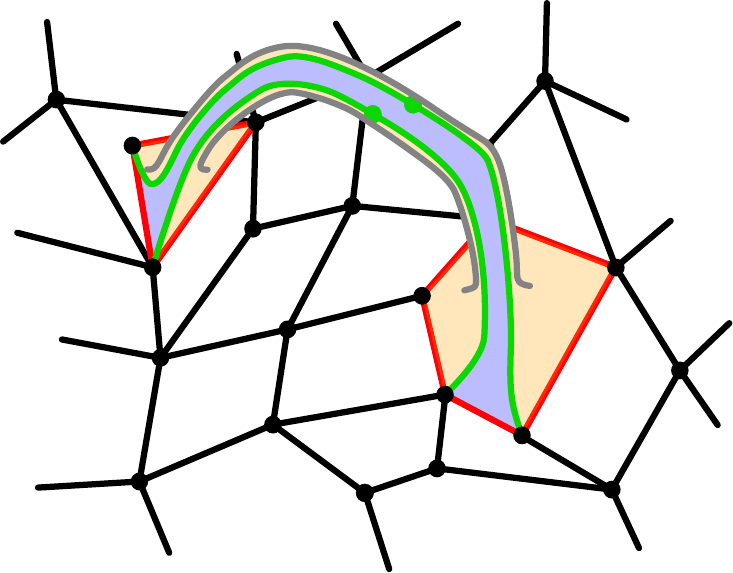}
\caption{An illustration of the process of adding a handle to connect a pair of odd faces. The left graph has exactly two odd faces, shown in red. In the right graph, the odd faces are connected by the gray handle, as well as the four green edges and two green vertices on the handle. Note that the resulting graph has two new faces, shaded in blue and orange, which are both even. Moreover, the red cycles do not bound faces anymore.}
\label{fig:oddfaces}
\end{figure}

We have now reduced to the case when $G$ is cellularly embedded in a surface of Euler genus $g$ with a constant number, say $k$, of odd faces. Clearly $k$ must be even by considering the number of odd degree vertices in the dual. We can then arbitrarily pair up these odd faces, and for each pair we add two vertex disjoint paths between them consisting of new vertices, and we give these new vertices weight 0, see \Cref{fig:oddfaces}. We think of these new paths as being embedded in a new handle added between the faces. We can do this such that the resulting graph has an even faced embedding in a surface of Euler genus $g + k$.

\subsection{Standard assumptions}\label{subsection:standard_assumptions}
Let $G$ be a graph with an even-faced embedding in a surface of Euler genus $g$. We now record a list of assumptions that we may assume when solving \Cref{problem:MIS_even_faced}, which we refer to as the ``standard assumptions''.

We say a vertex weighting $w: V(G) \rightarrow \mathbb{Q}_{\geq 0}$ is \emph{edge-induced} if there is a nonnegative edge weighting $c: E(G) \rightarrow \mathbb{Q}_{\geq 0}$ such that for every vertex $v$ we have $w(v) = \sum_{e \in \delta(v)} c(e)$.

\begin{assumption}\label{assumption:edge_induced_weights}
    The vertex weights are edge-induced, and we can compute the corresponding edge weighting in polynomial time.
\end{assumption}
\begin{proof}
    This follows from solving the linear programming relaxation, applying a result of Nemhauser and Trotter~\cite{nemhauser1974properties}, and complementary slackness, see \cite[Section 4]{conforti2020stable}. We note that in the reduction to edge-induced weights given by Conforti et al.~\cite{conforti2020stable} they possibly delete vertices, but we may instead set their weight to 0 (and set $c(e) = 0$ for all incident edges) in order to preserve a cellular embedding.
\end{proof}

\begin{assumption}
    $G$ has minimum degree 2.
\end{assumption}
\begin{proof}
    Clearly we may assume $G$ has no isolated vertices. Suppose $G$ contains a vertex $v$ of degree 1 with single neighbor $u$. Because weights are edge induced, $w(v) = c(uv) \geq 0$. Without loss of generality we can assume that one of $u$ or $v$ will be picked in a maximum weight independent set (even if $c(uv) = 0$). That is, $uv$ will never be slack. Then we may delete $v$ and add $c(uv)$ to our final weight (in doing so, the weight of $u$ decreases by $c(uv)$). By repeating this at most $|V(G)|$ times we reduce to an instance with minimum degree 2. We note that after deleting a vertex of degree 1 our embedding is still cellular and even-faced.
\end{proof}

\begin{assumption}
    $G$ is non-bipartite.
\end{assumption}
\begin{proof}
    It is well known that we can find a maximum weight independent set in $G$ in polynomial time if $G$ is bipartite.
\end{proof}

We note the above implies $G$ is not embedded in the plane.

\begin{assumption}
    $G$ is 2-connected and the embedding of $G$ has facewidth at least 2.
\end{assumption}
\begin{proof}

    We will use that for any cellularly embedded graph $H$, the faces of $H$ are all bounded by cycles if and only if $H$ is 2-connected with facewidth at least 2, see \cite[Proposition 3.8]{Mohar1997facewidth}. Hence, it suffices to force every face to be bounded by a cycle. For a vertex $v$ appearing on the boundary of a face $F$, we define the \emph{$F$-multiplicity} of $v$ to be the number of times that $v$ appears on the facial walk of $F$. We define the \emph{excess} of $v$ to be the sum over the faces $F$ containing $v$ on their boundary of the $F$-multiplicity of $v$ minus 1. We define the \emph{excess} of $G$ to be the sum of the excess of every vertex. If $G$ has excess 0, then every face is bounded by a cycle. We note that each vertex has excess at most its degree, so the excess of $G$ is bounded by $\mathcal{O}(|E(G)|)$.
    
    
    We do the following for each face $F$ of $G$. Consider a facial walk of $F$, and suppose that some vertex $v$ appears twice on this facial walk. Let $u$ and $w$ be the vertices preceding and following $v$ on this facial walk respectively. Because $G$ is simple and does not contain vertices of degree 1, $u \neq w$. We add a new vertex $v'$ of weight 0 adjacent to $u$ and $w$, and we embed $v'$ in the face $F$. By setting the cost $c$ of the edges $v'u$ and $v'w$ to be zero we preserve that the weights are edge induced, and clearly any optimal solution remains optimal. By adding the vertex $v'$ we have divided the even face $F$ into two even faces, and the resulting graph has smaller excess than $G$. Hence by repeating this a polynomial number of times we arrive at a graph of excess 0, and hence a graph where every face is bounded by a cycle. We note this preprocessing step is similar to one used in \cite[Section~6.4]{fiorini2025integer}.
\end{proof}


\subsection{Slack vectors}\label{subsection:slack_vectors}


As it turns out, in order to solve MIS instances on surfaces, it is useful to think of independent sets as integer vectors in the edge space. The goal of this section is characterize these vectors and explain how to formulate MIS in the edge space. Although most of the results follow from known results, in particular, from Conforti et al.~{\cite[Proposition~14, Lemma~16]{conforti2020stable}}, this section is self-contained, except for \Cref{lem:make_slack_01} which we do not reprove.

Let $G$ be a graph, and let $A \in\{0,1\}^{E(G)\times V(G)}$ denote its (edge-vertex) incidence matrix. Throughout this section, we assume that $G$ is connected, although the entire discussion easily generalizes to arbitrary graphs. Consider the set $\slackspace(G) :=  \{y \in \Q^{E(G)} \mid \exists x \in \Q^{V(G)} : Ax + y = \one\}$. This is an affine subspace of the edge space $\Q^{E(G)}$ that we call the \emph{slack space}. In fact, $\slackspace(G)$ is a linear subspace of the edge space, since $\zero \in \slackspace(G)$. Indeed, we have $A \one = 2 \cdot \one$ (the first $\one$ is the all-one vector in the vertex space and the second the all-one vector in the edge space). Hence, letting $x := \frac{1}{2} \one$ and $y := \zero$, we find $Ax + y =  A (\frac{1}{2} \one) + \zero = \one$, proving that $\zero \in \slackspace(G)$. This also shows that $\lambda^\intercal \one = 0$ for all $\lambda \in \ker(A^\intercal)$. (We point out that the above arguments break in characteristic $2$: it is important that we are working over $\Q$ and not say over $GF(2)$.) 

Our first lemma gives a characterization of $\slackspace(G)$ that follows directly from the Fredholm alternative.

\begin{lemma} \label{lem:Fredholm}
    The slack space $\slackspace(G)$ is the set of all $y \in \Q^{E(G)}$ such that $\lambda^\intercal y = 0$ for all $\lambda \in \ker(A^\intercal)$.
\end{lemma}


Next, we make \Cref{lem:Fredholm} more explicit and provide a non-redundant set of linear equations describing the slack space. An \emph{elementary vector} of $A^\intercal$ is a support-minimal element of $\ker(A^\intercal) \setminus \{\zero\}$. These are sometimes called \emph{minimal linear dependences} between the columns of $A^\intercal$. Note that, up to scaling, every elementary vector is determined by its support. The supports of elementary vectors coincide with the \emph{circuits} of the (column) vector matroid $M(A^\intercal)$. Every $\lambda \in \ker(A^\intercal)$ can be written as a linear combination of elementary vectors. In fact, we can even say more.

Let $B$ denote a (column) basis of $A^\intercal$. This is a maximal set of linearly independent columns of $A^\intercal$, or equivalently a basis of the vector matroid $M(A^\intercal)$. Every column of $A^\intercal$ that is not in the basis can be expressed in a unique way as a linear combination of the columns in the basis. The corresponding elementary vector $\lambda$ is called the \emph{fundamental elementary vectors} of $A^\intercal$ with respect to $B$. (By convention, we scale $\lambda$ so that $\lambda(j) = -1$, where $j$ is the label of the column that is not in the basis, and expressed in the basis.) Then, every element of $\ker(A^\intercal)$ can be written as a linear combination of the fundamental elementary vectors of $A^\intercal$ with respect to $B$. 

We conclude the following. Let $B$ be any basis of $A^\intercal$. Then, the slack space $\slackspace(G)$ is the set of all $y \in \Q^{E(G)}$ such that $\lambda^\intercal y = 0$ for all fundamental elementary vectors $\lambda$ of $A^\intercal$ with respect to $B$.

Now, since $A^\intercal$ is the vertex-edge incidence matrix of graph $G$, the bases and circuits of $M(A^\intercal)$, seen as subgraphs of $G$, are easy to describe. Since we assume $G$ is connected, the bases of $M(A^\intercal)$ are either the spanning trees (in case $G$ is bipartite), or a maximal set of edges in which every component contains at most 1 odd cycle (in case $G$ is non-bipartite). In the second case, every component of the basis is an odd $1$-tree. An \emph{odd $1$-tree} is a tree with a single extra edge whose unique cycle is odd. In particular, if $G$ is not bipartite, then it contains a basis formed from a spanning, odd 1-tree. The circuits of $M(A^\intercal)$ are exactly the even cycles and \emph{odd bicycles}, that is, the union of two odd cycles and a single path between them (the odd cycles are disjoint unless the path between them is trivial, in which case they share a single vertex). We describe the corresponding elementary vectors in terms of even closed walks.

A \emph{walk} in $G$ is a sequence $W=(v_0,e_1,v_1,e_2,v_2,\dots,e_k,v_k)$ where $v_0$, $v_1$, \ldots, $v_k$ are vertices and $e_1$, $e_2$, \ldots, $e_k$ are edges such that $e_i = v_{i-1}v_i$ for each $i \in [k]$. The \emph{length} of walk $W$ is $k$. Walk $W$ is \emph{closed} if $v_0 = v_k$. 
We call $W$ a \emph{cycle walk} if it is closed and $v_0$, \ldots, $v_{k-1}$ are distinct. 
We associate to walk $W$ the vector $\lambda_W := \sum_{i=1}^k(-1)^{i-1} \one_{e_i}$. Whenever $W$ is a closed even walk, we have $\lambda_W \in \ker(A^\intercal)$. Moreover, every elementary vector $\lambda$ of $A^\intercal$ is the form $\lambda = \lambda_W$, where $W$ is a walk around an even cycle or odd bicycle. We obtain the following non-redundant description of the slack space.

\begin{lemma} \label{lem:non-redundant_slack-space}
    Let $B$ denote any basis of $A^\intercal$. Let $\mathcal{W}$ denote the collection of even closed walks arising from fundamental circuits of $M(A^\intercal)$ with respect to $B$. Hence, $|\mathcal{W}| = |E(G)| - |V(G)| + 1$ if $G$ is bipartite, and $|\mathcal{W}| = |E(G)| - |V(G)|$ if $G$ is non-bipartite. The slack space $\slackspace(G)$ is the set of all $y \in \Q^{E(G)}$ such that $\lambda_W^\intercal y = 0$ for all $W \in \mathcal{W}$.
\end{lemma}

We call an even closed walk $W$ a \emph{circuit walk} if $W$ is a walk around an even cycle or an odd bicycle. Also, we define the \emph{width} of such a walk $W$ as $0$ in case $W$ is an even cycle walk, or as the number of vertices on the connecting path in case $W$ is a walk around an odd bicycle. 

\begin{lemma} \label{lem:circuit_walk_decomp}
    Let $W$ denote a circuit walk. If $G$ is $2$-connected, then $\lambda_W$ can be written as a linear combination of vectors $\lambda_{W'}$ where $W'$ is an even cycle walk.
\end{lemma}

\begin{proof}
    Let $k$ denote the width of $W$. If $k = 0$, then there is nothing to prove. Suppose that $k \ge 1$. 
    We show that $\lambda_W$ can be written as the sum of two vectors $\lambda_{W'}$ where $W'$ is a circuit walk of width at most $k-1$. The result follows by induction on $k$.

    Suppose we have an even closed walk $W$ around an odd bicycle consisting of two odd cycles $C_1, C_2$ and a path $P$ between them. Let $v_1$ and $v_2$ be the ends of $P$ in $C_1$ and $C_2$ respectively (possibly, $v_1 = v_2$, since $P$ may be trivial). Since $G$ is $2$-connected, we can find a path $Q$ from $C_1 - v_1$ to $(C_2 \cup P) - v_1$ that is internally disjoint from the odd bicycle. Let $w_1$ and $w_2$ denote the ends of $Q$ in $C_1 - v_1$ and $(C_2 \cup P) - v_1$, respectively. 

    If $w_2 \in V(C_2 - v_2)$, then we may find two even cycles $C'_1$ and $C'_2$ whose union is $C_1 \cup C_2 \cup P \cup Q$ and whose intersection is $P \cup Q$. Let $W'_1$ and $W'_2$ be corresponding even cycle walks. We can choose $W'_1$ and $W'_2$ such that since the signs of $\lambda_{W'_1}$ and $ \lambda_{W'_2}$ agree on $E(P)$ and disagree on $E(Q)$. Then $\lambda_W = \lambda_{W'_1} + \lambda_{W'_2}$.
    
    Otherwise, $w_2 \in V(P - v_1)$ (notice that this case can arise only when $k \ge 2$). Let $P'$ denote the subpath of $P$ from $v_1$ to $w_2$. Then we can find an even cycle $C'_1$ and an odd bicycle $C'_2$ such that the union of $C'_1$ and $C'_2$ is $C_1 \cup C_2 \cup P \cup Q$ and the intersection is $P' \cup Q$. Again, let $W'_1$ and $W'_2$ be corresponding circuit walks, thus $W'_1$ is an even cycle walk and $W'_2$ is an odd bicycle walk of width at most $k-1$. After possibly shifting one of $W'_1$ and $W'_2$, we can ensure that the signs of $\lambda_{W'_1}$ and $ \lambda_{W'_2}$ agree on $E(P')$ and disagree on $E(Q)$. Once more, we have $\lambda_W = \lambda_{W'_1} + \lambda_{W'_2}$.
\end{proof}

Together with \Cref{lem:non-redundant_slack-space},
\Cref{lem:circuit_walk_decomp} implies that the slack space has a complete description in the edge space via $|E(G)|-|V(G)|$ (resp.\ $|E(G)|-|V(G)|+1$) linear constraints that correspond to even cycles, whenever $G$ is $2$-connected and non-bipartite (resp.\ bipartite).

Next, consider the set $\slacklattice(G) := \{y \in \Z^{E(G)} \mid \exists x \in \Z^{V(G)} : Ax + y = \one\}$. This is an affine lattice that we call the \emph{slack lattice}. We call its elements \emph{slack vectors}. Observe that $\slacklattice(G)$ is contained in the slack space $\slackspace(G)$. 
Clearly, $\slacklattice(G) \subseteq \slackspace(G) \cap \Z^{E(G)}$. As stated in our next lemma, the converse inclusion holds whenever $G$ is bipartite, and when $G$ is non-bipartite adding a single parity constraint yields a complete description of $\slacklattice(G)$.

\begin{lemma} \label{lem:characterization_slack_lattice}
    If $G$ is bipartite, then $\slacklattice(G) = \Z^{E(G)} \cap \slackspace(G)$. If $G$ is non-bipartite, then $\slacklattice(G) = \{y \in \Z^{E(G)} \cap \slackspace(G) \mid \lambda_{W_0}^\intercal y \equiv 1 \pmod*{2}\}$, where $W_0$ is any walk around an odd cycle of $G$.
\end{lemma}

\begin{proof}
    First, assume that $G$ is bipartite, and let $y \in \Z^{E(G)} \cap \slackspace(G)$. Since $A$ is totally unimodular, there exists $x \in \Z^{V(G)}$ such that $Ax = \one - y$, showing that $y \in \slacklattice(G)$. 

    Second, assume that $G$ is non-bipartite. Fix some $y \in \Z^{E(G)} \cap \slackspace(G)$. We have to prove that $Ax = \one - y$ has a solution $x \in \Z^{V(G)}$ if and only if $\lambda_{W_0}^\intercal y \equiv 1 \pmod*{2}$. 
    
    Let $C_0$ denote the subgraph of $G$ traced by $W_0$. Let $B \subseteq G$ denote any odd, spanning $1$-tree containing $C_0$ (so $B$ is a basis). Consider the subsystem $A_B x = \one_B - y_B$, where the subscript $B$ indicates restriction to rows corresponding to edges of $B$. Since $G$ is non-bipartite, $A_B$ is invertible. We put $A_B$ in Hermite normal form. Assume that the elements of $B$ are ordered in such a way that the edges of the odd cycle come last, and in consecutive order. Then, there exists an unimodular matrix $U$ such that
    $$
    A_B U = 
        \begin{bmatrix}
        \identity & 0\\
        0 & C\\
        \end{bmatrix}
    \quad \text{where} \quad
    C = 
        \begin{bmatrix}
        1 &0 &\cdots &0 &0\\
        0 &1 & \cdots &0 &0\\        
        \vdots &\vdots& \ddots& \vdots &\vdots\\
        0 &0 &\cdots &1&0\\
        1 &1 &\cdots &1 &2
        \end{bmatrix}\,.   
    $$    
    Since $y \in \Z^{E(G)}$, we see that $A_B x = \one_B - y_B$ has a solution $x \in \Z^{V(G)}$ if and only if $\lambda_{W_0}^\intercal y \equiv 1 \pmod*{2}$. Since $y \in \slackspace(G)$, this occurs if and only if $A x = \one - y$ has an integer solution.

\end{proof}

Among other things, \Cref{lem:characterization_slack_lattice}
implies that $\zero \in \slacklattice(G)$ (that is, $\slacklattice(G)$ is a lattice) if and only $G$ is bipartite.


We wrap up the discussion of the slack lattice, slightly switching notation to be more in line with the literature. We recall the following definition from \Cref{section:algorithm_overview}: given walk $W=(v_0,e_1,v_1,e_2,v_2,\dots,e_k,v_k)$, we define a linear form $\omega_W$ on $\Q^{E(G)}$ by letting $\omega_W(y) := \sum_{i=1}^k(-1)^{i-1}y(e_i)$ for $y \in \Q^{E(G)}$. That is, we let $\omega_W(y) = \lambda_W^\intercal y$.

\begin{lemma}\label{lem:slacklattice-generators}
    For every $2$-connected graph $G$, there exists a collection $\mathcal{C}$ of $|E(G)| - |V(G)| + 1$ cycle walks such that     
    $$
    \slacklattice(G) = \{y \in \Z^{E(G)} \mid \begin{array}[t]{@{}l@{\ }l@{}}
    \omega_W(y) = 0 &\text{for every even cycle walk } W \in \mathcal{C},\\ 
    \omega_W(y) \equiv 1 \pmod{2} &\text{for every odd cycle walk } W \in \mathcal{C} \}\,.
    \end{array}
    $$
    If $G$ is non-bipartite, then $\mathcal{C}$ contains a unique odd cycle walk. Moreover, $\mathcal{C}$ can be constructed in polynomial time.
\end{lemma}

\begin{proof}
    The statement follows from our three previous lemmas. To see that the collection $\mathcal{C}$ can be computed in polynomial time, we start with a basis $\{\lambda_W \mid W \in \mathcal{W}(B)\}$ of $\ker(A^\intercal)$ where $\mathcal{W}(B)$ denotes the circuit walks coming from any odd, spanning $1$-tree $B$ (or spanning tree $B$, in case $G$ is non-bipartite). For each $W$ that is an odd bicycle walk, we perform at most $|E(G)|$ updates as in the proof of \Cref{lem:circuit_walk_decomp}. Each update finds two vectors $\lambda_{W'_1}, \lambda_{W'_2}$ which span $\lambda_{W}$. We remove $\lambda_W$ from the basis and replace it with one of $\lambda_{W'_1}, \lambda_{W'_2}$ such that we keep a basis. At the end of this process, we obtain a basis of $\ker(A^\intercal)$ of the form $\{\lambda_W \mid W \in \mathcal{C}\}$, where $\mathcal{C}$ is the resulting collection of even cycle walks. In case $G$ is non-bipartite, we pick any odd cycle walk $W_0$ and add it to $\mathcal{C}$.
\end{proof}

Finally, we discuss bounds on $x \in \Z^{V(G)}$ and $y \in \Z^{E(G)}$ on the solutions to $Ax + y = \one$. We require that $y \ge \zero$, which implies $Ax \le \one$. It would be natural to also enforce $x \ge \zero$, since only such integer points correspond to independent sets. However, we will not do this explicitly and rather show that this can be assumed without loss of generality whenever $(x,y)$ is optimal. We call a slack vector $y \in \slacklattice(G)$ \emph{feasible} if $y \ge \zero$.

Let $c \in \Q_{\ge 0}^{E(G)}$ and let $w \in \Q_{\ge 0}^{V(G)}$ satisfy $w(v) := c(\delta(v))$ for all $v \in V(G)$. This last condition can be written $w = A^\intercal c$. If $x \in \Z^{V(G)}$ and $y \in \Z^{E(G)}_{\ge 0}$ satisfy $Ax + y = \one$, then by multiplying on the left by $c^\intercal$ (and letting $w(x) := w^\intercal x$ and $c(y) := c^\intercal y$), we get $w(x) + c(y) = c(E(G))$. Notice that $c(E(G)) = \frac{1}{2}w(V(G))$. 

Our next lemma states that finding a maximum independent set in $G$ with respect to vertex weights $w(v)$ is equivalent to finding a minimum cost feasible slack vector with respect to edge costs $c(e)$. Notice that, by what precedes, we have $c(y') \le c(y)$ if and only if $w(x') \ge w(x)$ for solution pairs $(x,y)$ and $(x',y')$.

\begin{lemma}[Fiorini et al.~{\cite[Lemma~32]{fiorini2025integer}}] \label{lem:make_slack_01}
    Let $G$ be a connected, non-bipartite graph and let $w \in \Q_{\geq 0}^{V(G)}$ be vertex weights induced by edge costs $c \in \Q_{\geq 0}^{E(G)}$. Given any pair $(x,y) \in \Z^{V(G)} \times \Z_{\ge 0}^{E(G)}$ with $Ax + y = \one$, one can compute in strongly polynomial time a pair $(x',y') \in \{0,1\}^{V(G)} \times \{0,1\}^{E(G)}$ such that $Ax' + y' = \one$ and $c(y') \le c(y)$.
\end{lemma}

By the above we have reduced to \Cref{problem:min_cost_feasible_slack}, restated below.

\MinCostFeasibleSlack*


\section{The derived graph}\label{section:derived_graph}

An important object in our algorithm will be the \emph{derived graph} of a group-labeled graph. Derived graphs are important in topological graph theory. We provide a brief overview here for the specific case of signed graphs, but we refer the reader to \cite{gross2001topological} for a more in-depth introduction. The derived graph of a signed graph also appears in \cite{Zaslavsky1982SignedGraphs} where the relation between their matroid and the matroid of the signed graph is discussed.

Given a signed graph $(G,\gamma)$ where $\gamma: E(G) \rightarrow \mathbb{Z}_2$, we associate a \emph{derived (signed) graph} $(\tilde{G}, \tilde{\gamma})$ defined as follows. $\tilde{G}$ has vertex set $V(\tilde{G}) = V(G) \times \mathbb{Z}_2$. For each edge $uv \in E(G)$, we add the edges $(u,0)(v,\gamma(uv))$ and $(u,1)(v,1 + \gamma(uv))$ to $E(\tilde{G})$, and give these edges $\tilde{\gamma}$-label equal to $\gamma(uv)$. Intuitively, we have two copies of our vertex set, and by traversing an even edge one stays in the same copy and by traversing an odd edge one swaps to the other copy.

For each vertex $v \in V(G)$, we call the vertices $(v,0),(v,1) \in V(\tilde{G})$ its \emph{copies} in $\tilde{G}$, and we call $v$ the \emph{underlying vertex} of $(v,0)$ and $(v,1)$. Similarly for each edge $uv \in E(G)$, we call the edges $(u,0)(v,\gamma(uv))$ and $(u,1)(v,1 + \gamma(uv))$ in $E(\tilde{G})$ the \emph{copies} of $uv$, and we call the edge $uv$ the \emph{underlying edge} of $(u,0)(v,\gamma(uv))$ and $(u,1)(v,1 + \gamma(uv))$.

In the case of a signed graph where every edge has label $1 \in \mathbb{Z}_2$, the derived graph is also called a \emph{bipartite double cover}. We note that for any signed graph $(G,\gamma)$, the derived graph $(\tilde{G}, \tilde{\gamma})$ has no odd cycles (with respect to $\tilde{\gamma}$).

We will use that if $(G,\gamma)$ has an even-faced embedding on a surface of bounded genus, then $(\tilde{G},\tilde{\gamma})$ also embeds on a surface of bounded genus. Then we will find an alternating orientation $D^*$ of the dual $\tilde{G}^*$ and solve MIS on $G$ by solving a minimum cost circulation problem on $D^*$ with extra topological constraints.

Given an embedding scheme $(\pi,\lambda)$ of $(G,\gamma)$, we can define an embedding scheme $(\tilde{\pi},\tilde{\lambda})$ of $(\tilde{G},\tilde{\gamma})$, called the \emph{derived embedding}, in the following way. Let $v$ be a vertex of $\tilde{G}$ with incident edges $\tilde{e}_1, \dots, \tilde{e}_k$, and let $e_1, \dots, e_k$ be the corresponding underlying edges. Then we simply order the edges $\tilde{e}_1, \dots, \tilde{e}_k$ in $\tilde{\pi}$ to agree with the ordering of $e_1, \dots, e_k$ in $\pi$. We set $\tilde{\lambda}(\tilde{e}) = \lambda(e)$ for $e$ the underlying edge of $\tilde{e}$. We let $\tilde{\Sigma}$ be the surface defined by this derived embedding, and we call $\tilde{\Sigma}$ the \emph{derived surface}. We note that if $(G,\gamma)$ does not contain a $\gamma$-odd cycle, then $(\tilde{G}, \tilde{\gamma})$ will be two disjoint copies of $(G,\gamma)$, and $\tilde{\Sigma}$ will be disconnected.

\begin{theorem}[{\cite[Theorem 4.1.1]{gross2001topological}}]\label{thm:derived_graph_embedding}
    If $(G,\gamma)$ is a signed graph containing a $\gamma$-odd cycle with a cellular embedding in a surface $\Sigma$ of Euler genus $g$, then the derived graph $(\tilde{G}, \tilde{\gamma})$ embeds in a surface $\tilde{\Sigma}$ such that every even face of $(G,\gamma)$ becomes two faces of the same length in $(\tilde{G}, \tilde{\gamma})$, and every odd face of $(G,\gamma)$ becomes a single face in $(\tilde{G}, \tilde{\gamma})$ of twice the length. In particular, if $(G,\gamma)$ is connected, contains a $\gamma$-odd cycle, and has an even-faced embedding in $\Sigma$, then the Euler genus of $\tilde{\Sigma}$ is $2g-2$.
\end{theorem}

By the above it is well defined to talk of the two ``copies'' of an even face in $G$ as they appear in the derived embedding of $\tilde{G}$. We warn the reader that although every vertex of $G^*$ has two copies in the vertices of $(\tilde{G})^*$ when $G$ has an even-faced embedding, it is not true that $(\tilde{G})^*$ is the derived graph of $(G^*, \gamma)$. For instance if $(G,\mathbbm{1})$ is such that every edge is odd, then $\tilde{G}$ is bipartite, but its dual $(\tilde{G})^*$ may contain odd cycles whenever there is an odd degree vertex of $G$. Yet the derived graph of $(G^*,\mathbbm{1})$ is bipartite.

We also note the following.
\begin{theorem}[{\cite[Theorem 4.1.4]{gross2001topological}}]\label{thm:orientability_of_derived_embedding}
    Let $(G,\gamma)$ be a signed graph with embedding in $\Sigma$, and let $(\tilde{G}, \tilde{\gamma})$ be the derived graph with derived embedding in $\tilde{\Sigma}$. Then $\tilde{\Sigma}$ is non-orientable if and only if $(G,\gamma)$ contains a 1-sided even cycle.
\end{theorem}


We now describe how to obtain the surface $\tilde{\Sigma}$ from $\Sigma$ in terms of surface cutting in order to provide the reader with geometric intuition. We first note that if the signed graph $(G,\gamma')$ is obtained from $(G,\gamma)$ by shifting at the vertex $v \in V(G)$, then the derived graph of $(G,\gamma')$ is isomorphic to a shifting of the derived graph of $(G,\gamma)$. The derived graph of $(G,\gamma')$ can be obtained from the derived graph of $(G,\gamma)$ by exchanging the labels of $(v,0)$ and $(v,1)$ and shifting at those two vertices. Hence, up to shifting and isomorphism, the derived graph of a signed graph does not depend on the shifting.

\begin{figure}[ht]
\centering
\includegraphics[width=.45\textwidth]{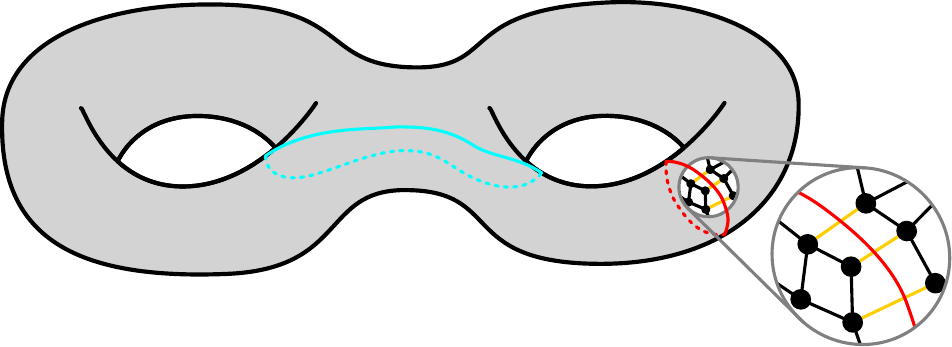}\hfill
\includegraphics[width=.45\textwidth]{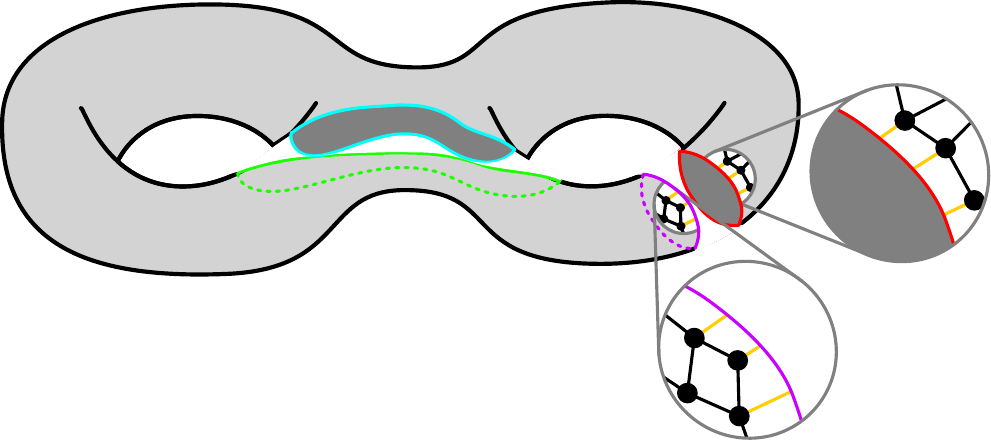}
\caption{An illustration of the process that lets us obtain the derived graph and the derived embedding starting from a graph $G$ embedded in an orientable surface of Euler genus $4$. The left figure shows the original embedded graph, where the odd edges are the yellow edges, i.e., the edges that cross the red and cyan curve. The right figure shows the process of cutting along the curves while maintaining the embedding of the graph.}
\label{fig:derived1}
\end{figure}

\begin{figure}[ht]
\centering
\includegraphics[width=.45\textwidth]{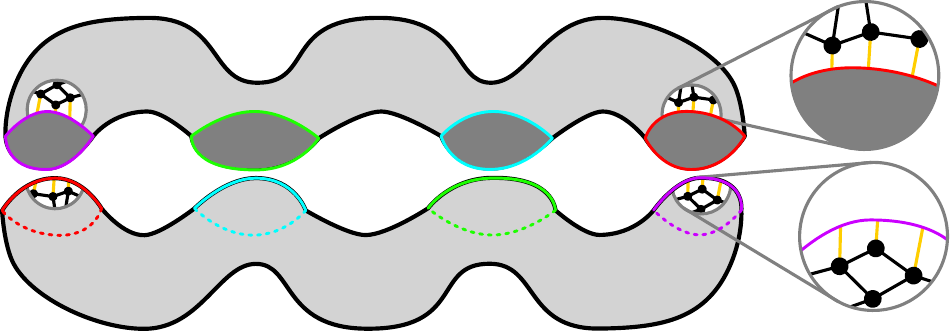}\hfill
\includegraphics[width=.45\textwidth]{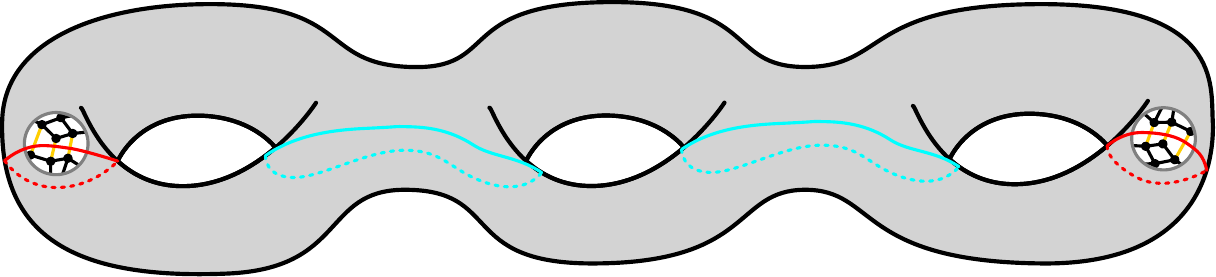}
\caption{In the left figure, we can see the two copies of the cut surface, aligned such that the corresponding embeddings fit together. The right figure shows the derived graph together with its derived embedding after gluing the copies at these curves.}
\label{fig:derived2}
\end{figure}

Given a signed graph $(G,\gamma)$ with an even-faced embedding in a surface $\Sigma$, the parity of a cycle depends only on the $\mathbb{Z}$-homology of the cycle in $\Sigma$. Thus we can find a set of closed curves in $\Sigma$ such that a cycle is odd if and only if it crosses these curves an odd number of times. Then we can shift $(G,\gamma)$ such that the odd edges are exactly the edges which cross these curves an odd number of times. The surface $\tilde{\Sigma}$ is formed from $\Sigma$ by cutting at these curves, copying the surface, and then gluing the copies at these curves. We illustrate the process in \Cref{fig:derived1} and \Cref{fig:derived2}.


\section{Reducing to a circulation problem}

Recall from \Cref{subsection:standard_assumptions} that the standard assumptions for a graph with associated weight function $(G,w)$ are the following:
\begin{itemize}
\item $G$ is embedded in a surface Euler genus $g = O(1)$ with all faces even,
\item $G$ is non-bipartite,
\item the embedding has facewidth at least $2$,
\item $G$ is $2$-connected, 
\item and the weight function $w$ is edge-induced.
\end{itemize}
The goal of this section is to reduce from \Cref{problem:min_cost_feasible_slack} to \Cref{problem:integer_circulation}, that is, to describe a feasible slack vector as a circulation over the dual of the derived graph.

Given an (unsigned, undirected) graph $G$ with an embedding scheme $(\pi,\lambda)$ in a surface $\Sigma$ of Euler genus $g$, we define the signed graph $(G, \mathbbm{1} + \lambda)$ where $(\mathbbm{1} + \lambda)(e) = 1 + \lambda(e) \in \mathbb{Z}_2$ for every edge $e \in E(G)$. We note that $(\pi, \lambda)$ also gives an even-faced embedding of $(G,\mathbbm{1} + \lambda)$ because every face of $G$ is homeomorphic to a disk, and hence the $\lambda$-weight of every facial cycle is 0. We denote its derived signed graph by $(\tilde{G}, \mathbbm{1} + \tilde{\lambda})$, and note that the label of every edge $e$ in the derived signed graph is given $1 + \tilde{\lambda}(e)$ where $(\tilde{\pi}, \tilde{\lambda})$ is the derived embedding in $\tilde{\Sigma}$.
By construction $(\tilde{G}, \mathbbm{1} + \tilde{\lambda})$ is balanced, i.e., for every cycle $C = v_1v_2\dots v_kv_1$ we have that $\sum_{i =1}^k (\mathbbm{1} + \tilde{\lambda})(v_iv_{i+1}) = 0$ (over $\mathbb{Z}_2$). This implies that the $\tilde{\lambda}$-weight of every cycle is equal to its parity (in the usual sense), or in other words, a cycle in $\tilde{G}$ is odd if and only if it is one-sided.

By \Cref{thm:derived_graph_embedding}, the derived embedding $(\tilde{\pi},\tilde{\lambda})$ is an embedding of $\tilde{G}$ on a surface $\tilde\Sigma$ of bounded genus. We take the dual $\tilde{G}^*$, and by \Cref{lem:alternating_orientation_equivalence} there exists an orientation $D^*$ of $\tilde{G}^*$ which is alternating. That is, for every vertex $v \in V(D^*)$, the edge incident to $v$ alternatively enter and leave $v$ in the cyclic ordering given by $\tilde{\pi}$.


It could be the case that $(G, \mathbbm{1} + \lambda)$ contains no odd cycles, that is, $G$ is embedded such that the odd cycles are exactly the one-sided cycles. We will call this case the \emph{degenerate case}. In this case, the derived graph $(\tilde{G}, \mathbbm{1} + \tilde{\lambda})$ is disconnected, and is formed from the disjoint union of two copies of $G$. We then work instead directly with the graph $G$ instead of the derived graph, and we define $D^*$ to be an alternating orientation of $G^*$. The proof becomes easier the degenerate case, and we will mention the small necessary changes throughout the proof. This is exactly the case which was solved previously by Conforti et al.~\cite{conforti2020stable}.


We will need the following property of the orientation $D^*$ when we are not in the degenerate case, and we note that the alternating orientation of $D^*$ is unique up to reversing every edge.

\begin{lemma}\label{lem:opposite_orientations_in_derived_dual}
    Suppose we are not in the degenerate case. Then for every edge $e \in E(G^*)$ with copies $e_0,e_1 \in E(\tilde{G}^*)$, the edges $e_0,e_1$ are given opposite orientations by $D^*$. That is, for $uv \in E(G^*)$ where $u_0,u_1$ are the copies of $u$ and $v_0,v_1$ are the copies of $v$ in $\tilde{G}^*$ such that $u_0v_0, u_1v_1 \in E(\tilde{G}^*)$, we have $\overrightarrow{u_0v_0},\overrightarrow{v_1u_1} \in E(D^*)$.
\end{lemma}
\begin{proof}
    Let $(\tilde{\pi},\tilde{\lambda})$ be the derived embedding of $\tilde{G}$, so in particular for every vertex $v \in V(G)$ its copies $v_0,v_1 \in V(\tilde{G})$ have the same cyclic ordering given by $\tilde{\pi}$. As in \Cref{lem:alternating_orientation_equivalence}, the orientation of $D^*$ is given by finding a shifting equivalent embedding scheme $(\tilde{\pi}', \tilde{\lambda}')$ such that $\tilde{\lambda}'(\tilde{e}) = 1$ for every edge $\tilde{e} \in E(\tilde{G})$, and then cyclically orienting the faces of $D^*$ to match with the cyclic orientation $\tilde{\pi}'$. Thus it suffices to show that for each vertex $v \in V(G)$ with copies $v_0,v_1 \in V(\tilde{G})$, the cyclic orientations given by $\tilde{\pi}'(v_0), \tilde{\pi}'(v_1)$ disagree. Recall that $V(\tilde{G}) = V(G) \times \mathbb{Z}_2$, and it follows from the definition of $\tilde{G}$ that $(\tilde{\pi}',\tilde{\lambda}')$ can be obtained from $(\tilde{\pi},\tilde{\lambda})$ by shifting at every vertex in $\{(v,1) \mid v \in V(G)\}$. Because the orientations given $\tilde{\pi}(v_0)$ and $\tilde{\pi}(v_1)$ agree and shifting swaps the orientation, the result follows.
\end{proof}

Recall that we have some edge weighting $c: E(G) \rightarrow \mathbb{Q}_{\geq 0}$. This naturally gives an edge weighting of $\tilde{G}$ (and hence of $D^*$) by giving each copy of an edge the same weight as the underlying edge. We define a linear map $\rho: \mathbb{Z}^{E(D^*)} \rightarrow \mathbb{Z}^{E(G)}$ such that entry $e \in E(G)$ of $\rho(f)$ has value $f(\vec{e}_0) + f(\vec{e}_1)$ for $e$ the underlying edge of both $e_0,e_1 \in E(\tilde{G})$. Then for any nonnegative integer circulation $f$ in $D^*$, we associate the vector $y \in \mathbb{Z}_{\geq 0}^{E(G)}$ given by $y = \rho(f)$. That is, $y(e) = f(\vec{e}_0) + f(\vec{e}_1)$ where $e_0,e_1 \in E(\tilde{G})$ are the two copies of $e \in E(G)$. Note that the orientations of $\vec{e}_0$ and $\vec{e}_1$ will disagree in $D^*$, but the flow does not cancel; both $ f(\vec{e}_0)$ and $f(\vec{e}_1)$ are nonnegative. In the degenerate case we define $\rho$ to be the identity function.

We will eventually solve a minimum cost nonnegative circulation problem on $D^*$ and return the associated vector $y = \rho(f)$. However, we first need to add some extra constraints by the following procedure. Recall from \Cref{subsection:slack_vectors} that $y \in \mathbb{Z}^{E(G)}_{\geq 0}$ is \emph{feasible} if $y \in \Lambda(G) = \{y \in \Z^{E(G)} \mid \exists x \in \Z^{V(G)} : Ax + y = \one\}$ where $A$ is the incidence matrix of $G$.

\begin{lemma}\label{lem:check_homology_generators}
    Let $G$ satisfy the standard assumptions and let $\mathcal{C}$ be the collection of cycles given in \Cref{problem:min_cost_feasible_slack}. Then we can in polynomial time find a set of cycles $C_1, \dots, C_g$ where $C_1$ is an odd cycle and $C_2, \dots, C_g$ are even cycles for $g$ the Euler genus of $\Sigma$ such that $y \in \mathbb{Z}_{\geq 0}^{E(G)}$ is feasible if and only if
    \begin{enumerate}
        \item $\omega_C(y) = 0$ for every facial cycle $C$,
        \item $\omega_{C_i}(y) = 0$ for every even cycle $C_2, \dots, C_g$, and
        \item $\omega_{C_1}(y) \equiv 1 \pmod 2$.
    \end{enumerate}
\end{lemma}
\begin{proof}
    Let $\mathcal{C} = \{C_1, \dots, C_t\}$ for $t = |E(G)| - |V(G)| + 1$ such that $C_1$ is the unique odd cycle in $\mathcal{C}$. Recall that $y \in \mathbb{Z}_{\geq 0}^{E(G)}$ is feasible if and only if $\omega_{C_1}(y) \equiv 1 \pmod 2$ and $\omega_{C_i}(y) = 0$ for all $C_2, \dots, C_t$. We note that $\{\omega_C(y) =0 \mid C\text{ a facial cycle}\}$ contains at least $|F(G)| - 1$ linearly independent constraints. To see this, note that any proper subset of the faces must contain some edge which is in only one face from the collection (because the dual is connected), so there can not be a linear dependence among less than $|F(G)|$ of the faces. In fact, all constraints in $\{\omega_C(y) =0 \mid C\text{ a facial cycle}\}$ are linearly independent unless we are in the degenerate case, see \Cref{lem:alternating_orientation_equivalence}.

    Thus we may extend these linearly independent constraints to a basis by including at most $(t-1) - (|F(G)| - 1) = -|V(G)| + |E(G)| - |F(G)| + 1 = g - 1$ of the constraints given by the cycles $C_2, \dots, C_{t}$ as desired. We note that except in the degenerate case, it suffices to take only a collection $C_1, \dots, C_{g-1}$. We will then add a linearly dependent constraint given by $C_g$ for notational convenience.

    
\end{proof}

Each constraint corresponding to a cycle $C_i$ above immediately gives a linear constraint for an integer circulation $f$ of $D^*$ because we set $y = \rho(f)$ as above. We then have the following problem.

\IntegerCirculation*


We now reduce \Cref{problem:min_cost_feasible_slack} to \Cref{problem:integer_circulation}.

\begin{theorem}\label{thm:reduction_to_integer_circulation}
    Let $f$ be a feasible solution to \Cref{problem:integer_circulation}. Then $\rho(f)$ is a feasible solution to \Cref{problem:min_cost_feasible_slack} with the same objective value. Conversely there exists an optimal solution $y\in\{0,1\}^{E(G)}$ to \Cref{problem:min_cost_feasible_slack} such that $y = \rho(f)$ for some feasible solution $f$ to \Cref{problem:integer_circulation}.
\end{theorem}
\begin{proof}
    Let $f$ be a feasible solution to \Cref{problem:integer_circulation}. We first show that $y := \rho(f)$ is a feasible solution to \Cref{problem:min_cost_feasible_slack}. We proceed by \Cref{lem:check_homology_generators}. Consider a facial cycle $C$ in $G$, and let $v_0, v_1$ be the associated vertices in $D^*$ that correspond to the two copies of $C$ in $\tilde{G}$. Let $e^1, e^2, \dots, e^\ell$ be the edges of $C$ in cyclic order, and let $e^i_0, e^i_1$ be the two copies of $e^i$ for $i \in [\ell]$. Then because $D^*$ is an alternating orientation and by flow conservation at $v_0$ and $v_1$ we have
    \begin{align*}
        0 &= \sum_{i=1}^\ell (-1)^{i-1}f(\vec{e}_0^i) + \sum_{i=1}^\ell (-1)^{i-1} f(\vec{e}_1^i)\\
        &= \sum_{i=1}^\ell (-1)^{i-1} y(e^i)\\
        &= \omega_C(y)
    \end{align*}
    as desired. We also have $\omega_{C_i}(y) = 0$ for every even cycle $C_i$ and $\omega_{C_1}(y) \equiv 1 \pmod 2$ by assumption. Thus by \Cref{lem:check_homology_generators}, $y$ is feasible, and clearly $y$ has the same objective value as $f$.

    Now let $y$ be an optimal solution to \Cref{problem:min_cost_feasible_slack}. By \Cref{lem:make_slack_01}, we may assume that $y \in \{0,1\}^{E(G)}$. We now construct a nonnegative integer circulation $f$ in $D^*$ such that $y = \rho(f)$. Let $F \subseteq E(G)$ be the set of edges for which $y$ takes value 1. Let $G^*$ be the dual of $G$, and for every vertex $v \in V(G^*)$ we locally assign $+1$ or $-1$ to the incident edges such that the assignment alternates as you traverse the edges clockwise around the vertex. This is possible because every vertex in $G^*$ has even degree. This assignment does not have to be consistent; it could be that $e$ has ends $u,v$ and $u$ assigns $e$ the value $+1$ and $v$ assigns the value $-1$. Let us denote the local signing of $e$ by $u$ as $\mathsf{sgn}_u(e) \in \{-1,1\}$. For $u \in V(G^*)$, let $C$ be the facial walk around the face $u$ in $G$. Then $\omega_C(y) = 0$ implies $\sum_{e \in F \cap \delta(u)} \mathsf{sgn}_u(e) = 0$. We will now construct a collection of closed walks in $G^*$ using only edges in $F$ such that each edge of $F$ is used exactly once over all walks in the collection.

    We start at an arbitrary vertex $v$ incident to an edge in $F$, and walk along an arbitrary edge of $F$. Each time we enter a vertex $u$, we look at the local signing of our edge at $u$ and exit through an edge of the opposite sign. If $u \neq v$, then this is always possible because $\sum_{e \in F \cap \delta(u)}\mathsf{sgn}_u(e) = 0$. We continue like this for as long as possible. In particular, we continue even if we again reach the starting vertex $v$ as long as there is an edge of the opposite sign from which to leave $v$. If we ever get stuck, then it must be because we reached the starting vertex $v$, and furthermore it must be that the last edge of the walk has the opposite sign from the first edge of our walk. Otherwise, suppose without loss of generality the first traversed edge $e$ had $\mathsf{sgn}_v(e) = 1$, then the edges of $F \cap \delta(v)$ traversed so far have alternating sum equal to 2, so there is some untraversed edge $e'$ in $F \cap \delta(v)$ with $\mathsf{sgn}_v(e') = -1$ for which we may continue the walk.

    We may not have used all of the edges of $F$ in the above closed walk. However we can repeat the above as many times as needed, only using new edges of $F$ in each walk, to obtain the desired collection of closed walks.

    We note that for a vertex $u$, the edges of $\delta(u)$ traversed in our walk do not have to have strictly alternating signs with respect to $\mathsf{sgn}_u$. We only require that when we enter through an edge of one sign, we exit through an edge of the opposite sign. It could be, for instance, that the edges of $\delta(u)$ are traversed in the order $+, -,-,+$. However, it will be the case that whenever the alternation switches, this will correspond to walking to the other copy of $u$ in $D^*$, which we now describe.

    We now claim that each walk in our collection corresponds naturally to a closed walk in $D^*$ which obeys the direction of edges. Then by pushing a unit of flow along each walk in our collection we obtain a nonnegative integer circulation $f$ of $D^*$ such that $\rho(f) = y$.
    Let $v$ be the first vertex of our walk with copies $v_0,v_1$ in $D^*$, and by possibly reversing our walk we may assume that the first edge of our walk is oriented out of $v_0$, and we can start our walk there. Then if we enter a vertex $u$ through a directed edge, the next edge of our walk will correspond to an edge directed out of $u$ because the orientation of $D^*$ is alternating. Thus the walk in $G^*$ corresponds to a walk respecting directions in $D^*$, and it suffices to show that the walk in $D^*$ is closed. This is clear in the degenerate case, so suppose we are not in that case. Then by \Cref{lem:opposite_orientations_in_derived_dual} recall that for every edge $e \in E(G)$, its two copies $e_0,e_1$ have different orientations in $D^*$. Suppose towards contradiction that we end our walk at the copy $v_1$ of $v$ instead of $v_0$. Then the first edge of our walk is oriented out of $v_0$ and the last edge is oriented into $v_1$. However because the edges incident to $v_0$ and $v_1$ have opposite orientations, this implies that the first and last edge of our walk have the same local signing relative to $v$, a contradiction.
    

    Because $y$ satisfies $\omega_{C_i}(y) = 0$ for every even cycle $C_i$ and $\omega_{C_1}(y) \equiv 1 \pmod 2$, we have that $\rho(f) = y$ implies that $f$ is a feasible solution to \Cref{problem:integer_circulation}.
\end{proof}

\section{Final algorithm}\label{section:final_algorithm}

The main result of this section is that we may solve \Cref{problem:integer_circulation} in polynomial time. We will solve \Cref{problem:integer_circulation} via the following algorithm, which is due to Chambers, Erickson, and Nayyeri~\cite{chambers2012homology} in the orientable case and due to Morell, Seidel, and Weltge~\cite{morell2021minimum} in the nonorientable case.

\begin{theorem}[{\cite{chambers2012homology, morell2021minimum}}]\label{thm:min_cost_circulation_of_fixed_homology}
    There exists a polynomial time algorithm which takes as input a directed graph $D$ cellularly embedded on a surface of bounded genus with nonnegative costs $c: E(D) \rightarrow \mathbb{Q}_{\geq 0}$ and an integer circulation $f$ in $D$, and finds a minimum-cost nonnegative integer circulation in $D$ that is $\mathbb{Z}$-homologous to $f$.
\end{theorem}

\begin{theorem}
    \Cref{problem:integer_circulation} can be solved in polynomial time.
\end{theorem}

\begin{proof}
    For $f$ a nonnegative integer circulation in $D^*$, we set $\omega(f) = (\omega_{C_1}(\rho(f)), \omega_{C_2}(\rho(f)), \dots, \omega_{C_g}(\rho(f))) \in \mathbb{Z}_2 \times \mathbb{Z}^{g-1}$. We first claim that $\omega(f)$ depends only on the $\mathbb{Z}$-homology of $f$. Because $\omega$ is linear in $f$, it suffices to show that if $g$ is a circulation consisting of a single unit of flow around a single face, then $\omega(g) = \mathbf{0}$. Let $g$ be a circulation sending a single unit of flow around a face $C$ of $D^*$ corresponding to vertex $\tilde{v} \in V(\tilde{G})$. Let $v$ be the underlying vertex in $G$, so $\rho(g)$ is the indicator of $\delta(v)$. Then for each $i \in [g] \setminus \{1\}$, the cycle $C_i$ either doesn't include $v$ or includes two edges of $\delta(v)$ which appear with opposite signs in $\omega_{C_i}$. In either case we see that $\omega_{C_i}(\rho(g)) = 0$. Because $C_1$ either doesn't include $v$ or includes two edges in $\delta(v)$, we also have $\omega_{C_1}(\rho(g)) = 0 \in \mathbb{Z}_2$ as desired.

    Because the feasibility of $f$ depends only on its $\mathbb{Z}$-homology, we may guess the correct $\mathbb{Z}$-homology of $f$ and compute a minimum cost circulation with that $\mathbb{Z}$-homology by \Cref{thm:min_cost_circulation_of_fixed_homology}. The main difficulties are obtaining a polynomially bounded region over which to guess, and constructing an integer flow with the guessed $\mathbb{Z}$-homology in order to apply \Cref{thm:min_cost_circulation_of_fixed_homology}. Let $\tilde{g}$ denote the Euler genus of $\tilde{\Sigma}$.
    
    By~\cite[Corollary~12, Proposition~18]{morell2021minimum}, we can efficiently compute a collection of vectors $\xi_1, \dots, \xi_{\tilde{g}} \in \mathbb{Z}^{E(D^*)}$ such that the entries of $\xi_i$ are bounded by 2 in absolute value for every $i \in [\tilde{g}]$, and two circulations $f,f'$ of $D^*$ are homologous if and only if 
    \begin{itemize}
        \item $\langle \xi_i, f\rangle = \langle \xi_i, f'\rangle$ for each $i \in [\tilde{g}]$ if $\tilde{\Sigma}$ is orientable, or
        \item $\langle \xi_1, f\rangle \equiv \langle \xi_1, f'\rangle \pmod 2$ and $\langle \xi_i, f\rangle = \langle \xi_i, f'\rangle$ for all $i \in [\tilde{g}]\setminus \{1\}$ if $\tilde{\Sigma}$ is non-orientable.
    \end{itemize}
    We then define $\xi(f) = (\langle \xi_1, f\rangle, \dots, \langle \xi_{\tilde{g}}, f\rangle)$ either as a vector in $\mathbb{Z}^{\tilde{g}}$ if $\tilde{\Sigma}$ is orientable, or as a vector in $\mathbb{Z}_2\times \mathbb{Z}^{\tilde{g} - 1}$ if $\tilde{\Sigma}$ is non-orientable. That is, $\xi$ defines an isomorphism from the first homology group of $\tilde{\Sigma}$ over the integers to $\mathbb{Z}^{\tilde{g}}$ or $\mathbb{Z}_2\times \mathbb{Z}^{\tilde{g} - 1}$.
    
    By \Cref{thm:reduction_to_integer_circulation}, there is an optimal solution $f^*$ to \Cref{problem:integer_circulation} such that $\rho(f^*)=y$ for $y \in \{0,1\}^{E(G)}$. By definition of $\rho$, we obtain that there exists an optimal solution $f^* \in \{0,1\}^{E(D^*)}$. Because $\xi_i$ has entries bounded by 2 in absolute value for every $i$, we have that $\xi(f^*)$ has entries bounded by $2|E(D^*)|$ in absolute value (where if the first entry is in $\mathbb{Z}_2$, there is no restriction on this entry). Let $\Omega$ denote the set of all such $\xi(f^*)$, so $|\Omega| \leq (4|E(D^*)| + 1)^{\tilde{g}}$.

    It now suffices, for each $z \in \Omega$, to compute a circulation $f_z \in \mathbb{Z}^{E(D^*)}$ with $\xi(f_z) = z$. Suppose we can find such an $f_z$ for each $z$. Then we check whether $\omega(f_z) = (1,0,\dots,0)$. If so, we apply \Cref{thm:min_cost_circulation_of_fixed_homology} to find a minimum cost nonnegative integer circulation homologous to $f_z$, and we output the circulation of minimum cost over all such $z$.

    We now show how to compute an integer circulation $f_z$ with $\xi(f_z) = z$ for each $z \in \Omega$. We first compute a collection of cycles in $\tilde{G}$ which generate the first homology group of $\tilde{\Sigma}$ over the integers.
    This can be done by taking a spanning tree $T^*$ of $\tilde{G}^*$ and a spanning tree $T$ of $\tilde{G}$ which are edge disjoint; this is always possible and can be done in polynomial time \cite[Lemma 3.1]{eppstein2003dynamic}. Then the fundamental cycles of edges $e \in E(\tilde{G}^*) \setminus (E(T) \cup E(T^*))$ with respect to $T^*$ generate the fundamental group of $\tilde{\Sigma}$ \cite[Lemma 3.2]{eppstein2003dynamic}. Because the first homology group of $\tilde{\Sigma}$ over the integers is the abelianization of its fundamental group, these cycles also generate the first homology group of $\tilde{\Sigma}$ over the integers. That is, by sending one unit of flow around these cycles in $D^*$ (so that edges in the reverse direction of our cyclic orientation get negative flow), we have that any circulation is an integer combination of these cycles and facial circulations. Note that we have $|E(\tilde{G}^*)| - (|V(\tilde{G}^*)| - 1) - (|F(G^*)| - 1) = \tilde{g}$ such cycles, and let us denote them by $B_1, \dots, B_{\tilde{g}} \in \mathbb{Z}^{E(D^*)}$. 
    
    Because $B_1, \dots, B_{\tilde{g}}$ generate the first homology group of $\tilde{\Sigma}$, for each $z \in \Omega$ there exists an integer combination $a \in \mathbb{Z}^{\tilde{g}}$ such that $\xi(\sum_{i =1}^{\tilde{g}} a_iB_i) = z$. We can solve for $a$ either by finding an integer solution to the linear system
    \begin{align*}
        \sum_{i=1}^{\tilde{g}} a_i \xi(B_i) = z
    \end{align*}
    in the case that $\tilde{\Sigma}$ is orientable, or by finding an integer solution $(a_1, \dots, a_{\tilde{g}},k)$ to the linear system
    \begin{align*}
        \sum_{i = 1}^{\tilde{g}} a_i\langle \xi_1, B_i\rangle &= z_1 + 2k\\
        \sum_{i = 1}^{\tilde{g}} a_i\langle \xi_2, B_i\rangle &= z_2\\
        &\vdots\\
        \sum_{i = 1}^{\tilde{g}} a_i\langle \xi_{\tilde{g}}, B_i\rangle &= z_{\tilde{g}}\\
    \end{align*}
    where $z_1 \in \mathbb{Z}_2$ is interpreted as an element of $\{0,1\}$ in the case that $\tilde{\Sigma}$ is non-orientable. We may then set $f_z = \sum_{i =1}^{\tilde{g}} a_iB_i$ as desired.
\end{proof}
In the above proof we use \Cref{thm:min_cost_circulation_of_fixed_homology} as a black box. It's not hard to see that in \cite[Section 3]{morell2021minimum}, they prove \Cref{thm:min_cost_circulation_of_fixed_homology} by solving a minimum cost circulation with a constant number of extra linear constraints where one of these constraints is taken modulo 2. That is, the same proof as in \cite[Section 3]{morell2021minimum} can be used to solve \Cref{problem:integer_circulation} directly.

\phantom{}\\
\noindent\textbf{AI disclosure.}
No AI tools were used in the writing of this manuscript. All proofs and their presentations are entirely due to the authors.

\phantom{}\\
\noindent
\textbf{Acknowledgments.}
We thank Mujin Choi, Maximilian Gorsky, Gunwoo Kim, and Sebastian Wiederrecht for valuable discussions.

\bibliographystyle{siamplain}
\bibliography{references}

\end{document}